\documentclass[12pt,oneside]{article}
\newcommand{\Path }{}
\usepackage{\Path Paper_style} % Loads my formatting 
\usepackage{\Path Keywords_WB} % Loads my keywords
\usepackage{\Path Environments} % Loads my keywords
\usepackage{sidecap}
\usepackage{tikz-cd} 
\usepackage{natbib}
\usepackage{\Path tikzcd2}
\newcommand{\Fol}{\mathcal{F}}

\newcommand{\MesL}{\nu}
\newcommand{\Rect}{\mathcal{R}}

\newcommand{\Return}{R}

\newcommand{\Iso}{\Psi}
\newcommand{\IsoInv}{\Xi}

\newcommand{\Susp}{\Gamma}
\newcommand{\metric}{\tau}

\begin{document}
\title{\vspace{-1cm} Smooth Koopman eigenfunctions}
\author{Suddhasattwa Das\footnotemark[1]}
\footnotetext[1]{Department of Mathematics and Statistics, Texas Tech University, Texas, USA}
\date{\today}
\maketitle
\begin{abstract} 
	Any dynamical system, whether it is generated by a differential equation or a transformation map on a manifold, induces a dynamics on functional-spaces. The choice of functional-space may vary, but the induced dynamics is always linear, and codified by the Koopman operator. The eigenfunctions of the Koopman operator are of extreme importance in the study of the dynamics. They provide a clear distinction between the mixing and non-mixing components of the dynamics, and also reveal embedded toral rotations. The usual choice of functional-space is $L^2$, a class of square integrable functions. A fundamental problem with eigenfunctions in $L^2$ is that they are often extremely discontinuous, particularly if the system is chaotic. There are some prototypical systems called skew-product dynamics in which $L^2$ Koopman eigenfunctions are also smooth. The article shows that under general assumptions on an ergodic system, these prototypical examples are the only possibility. Moreover, the smooth eigenfunctions can be used to create a change of variables which explicitly characterizes the weakly mixing component too.
\end{abstract}

\begin{keywords}
Koopman operator, skew-product, quasiperiodicity, Discrete spectrum, submersion
\end{keywords}
\begin{AMS}	37A05, 37A10, 37A15, 37A25, 37A45, 37C10, 37C55, 37C85, \end{AMS}

%-_-_-_-_-_-_-_-_-_-_-_-_-_-_-_-_-_-_-_-_-_-_-_-_-_-_-_-_-_-_-_-_-_-_-_-_-_-_-_-_-_-_-_-_-_-_-_-_-_-_-_-_-_-_-_-_-_-_-_-_-_-_-_-_-_-_-_-_-_-_-_-_-_-_-_-
\section{Introduction} \label{sec:intro}

Dynamical systems theory involves the study of a map $\Phi : M\to M$ on a topological space $M$ or a continuous family of maps $\Phi^t:M\to M$. Although this branch of mathematics had started out as a descriptive tool for various physical phenomenon, the realization that the map(s) can induce a variety of different mathematical behavior within $M$ lead to a great diversification of the theory. The following will be our standing assumption : 

\begin{Assumption}\label{A:A1}
	$M$ is an $n$-dimensional $C^1$ manifold equipped with the the Borel $\sigma$-algebra $\Borel$. $\Phi^t$ is a $C^1$ flow on $M$ generated by a vector field $V$. The flow has an invariant probability measure $\mu$ with a compact support $X\subset M$.
\end{Assumption}

There has been mainly three means of studying dynamics - as a topological system \citep[][e.g.]{AkinEtAl1996, Blender2, DasJim17_chaos}, as a measure preserving system  \citep[][e.g.]{Halmos1944mixing, Halmos1944, Halmos1956, Walters2000}, or indirectly as a group of linear transformations induced on some function space \citep[][e.g.]{Nadkarni, DasGiannakis_RKHS_2018}. In all these approaches, one could encounter different facets of the phenomenon of \emph{chaos} \citep[][e.g.]{DGJ_compactV_2018, Katok_periodic_1980, Bowen_AxiomAflow_1975, Bowen_equi_1975, Robinson1998, Bowen_entropy_1970}. In a series of independent works spanning decades, it has been discovered that chaotic systems have a dense collection of unstable hyperbolic periodic orbits \citep[][e.g.]{Katok_periodic_1980}, and these points typically have dense stable and unstable foliation \citep[][e.g.]{Pesin_families_1976}. This results in a non-smooth geometry of the attractor which is locally dense with expanding and contracting regions, and in spite of this great complexity, can be approximated well by arbitrarily long periodic orbits \citep[][e.g.]{WangSun2010, LiaoEtAl2018} embedded within the dynamics.

There is yet another broad category of dynamical systems which are completely free of these traits of chaos. These are called \emph{quasiperiodic systems} \citep[][e.g.]{Dioph_Herman_1979, DasJim2017_SuperC, DSSY2017_QQ}. They are characterized by zero Lyapunov exponents, absence of stable or unstable manifolds, lack of mixing, and every single trajectories being dense. Topologically they are conjugate to a rotation on a torus. One of the long standing conjectures in dynamics  \cite{SanderYorkeChaos2015} is that a typical non-periodic, minimal dynamical system is either quasiperiodic or chaotic set. The operator theoretic formulation of dynamical systems provides a clearer distinction between these two types of dynamics.

\paragraph{The Koopman operator.} Koopman operators \cite{Nadkarni} act on observables by composition with the flow map. Given a function $f:M\to \cmplx$ or $f:X\to\cmplx$, and a time $t\in\real$, the action of the Koopman operator $U^t$ is defined as the time shifted map
\begin{equation}\label{eqn:def_Koop}
(U^t f)(x) := f\left( \Phi^t x \right).
\end{equation}
This makes $U^t$ a linear operator on the vector spaces of functions on $X$ or $M$. This action can be factored into the space of Borel functions on $X$ or $M$. $U^t$ also preserves the Banach spaces $C(X)$ and $C^r_c(M)$, where $C^r_c(M)$ is the space of compactly supported $C^r$ functions on $M$. $U^t$ also preserves the Banach spaces $L^p(\mu)$ for $p>0$ and the invariant measure $\mu$. Moreover, $U^t$ is a norm-preserving operator on these Banach spaces. Let $V$ be the vector field generating the flow. $V$ acts as a differential operator $V:C^1(M) \to C(M)$ , defined as
\begin{equation*}\label{eqn:def_gen_flow}
V f:=\lim_{t\to 0} \frac{ 1 }{ t } \left(U^t f - f\right), \quad f \in C^1(M).
\end{equation*}
When $p=2$, $L^2(\mu)$ is a Hilbert space, and $\{ U^t : t\in\real \}$ is a 1-parameter unitary evolution group on $L^2(\mu)$ \cite{DasGiannakis_delay_2019, DGJ_compactV_2018}. Such a group has a generator $\hat{V}$ acting on some dense subspace of $L^2(\mu)$, and the action of $\hat{V}$ on $C^1(M)\cap L^2(\mu)$ coincides with that of $V$. This space $L^2(\mu)$ will be our primary choice of functional space, in which we study the spectrum of the Koopman operator. The Hilbert space property provides a powerful array of tools to model the Koopman operator, both theoretically and numerically \cite[e.g.]{FroylandEtAl14, FroylandEtAl_coherent_2010, DasGiannakis_delay_2019, DGJ_compactV_2018, Philipp2024error}. An important property to consider for both the Koopman group and its generator, is the presence of eigenfunctions.

\paragraph{Koopman eigenfunctions} Since $U^t$  is a linear operator on $C^r(M)$, $C(X)$ and $L^2(\mu)$, one can  look for Koopman eigenfunctions in either of these spaces. Based on \eqref{eqn:def_Koop}, an eigenfunction $z:M\to\cmplx$ of $U^t$,  with eigenvalue $\lambda$ should satisfy an equation of the form
\begin{equation}\label{eqn:Def_koop_eigen}
U^t z (x) = z(\Phi^t x) = \lambda z(x) 
\end{equation}
The function $z$ above will be called a $C^r(M)$-eigenfunction if $z\in C^r(M)$ and \eqref{eqn:Def_koop_eigen} holds for every $x\in M$. Similarly, $z$ will be called a $L^2(\mu)$-eigenfunction if $z\in L^2(\mu)$ and \eqref{eqn:Def_koop_eigen} holds in $L^2(\mu)$ sense, i.e., $\norm{ U^t z - \lambda z}_{L^2(\mu)} = 0$. In this Hilbert space setting, every eigenfunction has associated to it a quantity $\omega$ called its \emph{eigenfrequency}, which is independent of the time $t$, and which satisfies
\begin{equation}\label{eqn:Koop_eigen_L2}
U^t z (x) = z(\Phi^t x) = e^{\iota \omega t} z(x) , \quad \forall t\in\real, \quad \mu-\text{a.e.} \ \ x\in M
\end{equation}
This is consistent with the fact that $U^t$ is a unitary operator on $L^2(\mu)$ and all its eigenvalues must have norm $1$. Now note that if $z$ is a $C(M)$-eigenfunction, then it is an $L^2(\mu)$ eigenfunction too for every finite invariant measure $\mu$, thus \eqref{eqn:Def_koop_eigen} would still hold with $\lambda=e^{\iota \omega t}$.  An important link between the unitary group $U^t$ and its generator $V$ is that they share the same $L^2(\mu)$ eigenfunctions :
\[ U^tz=\exp(\iota \omega t)z \quad \Leftrightarrow \quad Vz=\iota \omega z. \]
Koopman eigenfunctions carry significance in several real world phenomenon, such as latent seasonal sources in the study of timeseries \cite{DasMustAgar2023_qpd, DasEtAl2023traffic}, as indicators of coherence in fluid flows \cite{FroylandEtAl14b, DellnitzEtAl00}, and even as states in the operator theoretic description of quantum mechanics \citep[e.g.][]{Eisner2017pointwise}, and in the study of general evolution semigroups \cite{Suchanecki2010evolution, DGJ_compactV_2018, EisnerEtAl15}. Note that the constant functions are always eigenfunctions, corresponding to eigenvalue 1 (or eigenfrequency 0). Non-trivial eigenfunctions may or may not be present. If they happen to be present, they provide a meaningful splitting of the function space.

\paragraph{Spectral splitting} Let $\Disc$ denote the subspace of $L^2(\mu)$ spanned by the eigenfunctions of $V$ with non-zero eigenvalues, $\Disc^\bot$ denote its orthogonal complement. 
\begin{equation} \label{eqn:L2_decomp}
L^2(\mu)=\Disc\oplus\Disc^\bot.
\end{equation}
Note that $\Disc$ always contains the constant functions, which correspond to eigenfrequency $0$. If $d\geq 1$, there is at least one non-trivial Koopman eigenfunction. The eigenvalues of $V$ are closed under integer linear combinations and is generated by a finite set of rationally independent eigenvalues $\iota \omega_1,\ldots,\iota \omega_d$. Thus every eigenvalue of $V$ is of the form $\Sigma_{j\in\braces{ 1, \ldots, d}}a_j \omega_j$, for some $(a_1,\ldots,a_d)\in\integer^d$. If $d>1$, then the eigenvalues are dense on the imaginary axis. For every $\vec{a}=(a_1,\ldots,a_d)\in\integer^d$, $z_{\vec{a}}$ will denote the eigenfunction $\Sigma_{j\in\braces{ 1, \ldots, d}}a_j \omega_j$. The corresponding eigenvalue is $i\vec{a}\cdot\vec{\omega}$, where $\vec{\omega}$=$(\omega_1,\ldots,\omega_d)$.

The splitting in \eqref{eqn:L2_decomp} is yet another manifestation of the dichotomy of chaotic and quasiperiodic dynamics. By a reuse of terminology, the subspaces $\Disc$ and $\Disc^\bot$ will be called the quasiperiodic and chaotic components of the dynamics respectively. The action of $U^t$ on the quasiperiodic component is characterized by absence of mixing, unique ergodicity and superior ergodic convergence properties \cite[e.g.]{DasJim2017_SuperC,DSSY2017_QQ,DasYorkeQuasiR2016}. Chaotic motion is opposite in nature to quasiperiodicity \cite[e.g.]{DasJim17_chaos,DasToungue2018,Blender2}, it is characterized by infinitely many invariant subsets, and ergodic properties like arbitrarily small rate of decay of Birkhoff averages, mixing, and decay of correlations. We next further explore the connection of $\Disc$ with toral rotations. 

\paragraph{Embedded toral rotations.} In what follows, the circle $S^1$ is identified with the unit circle $\{ e^{i\theta} \; : \; \theta\in [0,2\pi)\}$ in the complex plane. If $z:M\to\cmplx$ is a Koopman eigenfunction nonzero everywhere, then $z/|z|$ is a Koopman eigenfunction too, with the same eigenfrequency. Thus without loss of generality (WLOG), $z$ can be considered to map into the circle, i.e., $z:M\to S^1$. For every $d\in\num$, the $d$-dimensional torus $\TorusD{d}$ is the $d$-times Cartesian product  $S^1 \times \cdots \times S^1$. An important property of Koopman eigenfunctions is that they factor the dynamics onto a rotation on $S^1$ with frequency $\omega$. This is shown on the left below : 
\begin{equation}\label{eqn:pi_factor}
	\begin{tikzcd}[row sep = small]
		M \arrow[swap]{d}{z\ } \arrow{r}{\Phi^t_V} &M \arrow{d}{\ z} \\
		S^1 \arrow{r}{R_\omega^t} &S^1
	\end{tikzcd}; \quad
	R_\omega^t(\theta) \mapsto \theta+t\omega \bmod S^1; \quad
	\begin{tikzcd}[row sep = small]
		M \arrow[swap]{d}{(z_1,\ldots,z_d)\ } \arrow{r}{\Phi^t_V} &M \arrow{d}{\ (z_1,\ldots,z_d)} \\
		\TorusD{d} \arrow{r}{R_\omega^t} &\TorusD{d}
	\end{tikzcd}; \quad
	R_\omega^t(\theta) \mapsto \theta+t\omega \bmod \TorusD{d}; \quad
\end{equation}
Similarly, $d$ Koopman eigenfunctions $(z_1,\ldots,z_d)$ factor the dynamics into a rotation on $\TorusD{d}$, with rotation vector $\left( \omega_1,\ldots,\omega_d \right)$. This is shown on the right above. If $z$ was an $L^2(\mu)$ eigenfunction instead of $C^1(M)$, then the factorization would be of $X$, and would hold in a measure theoretic sense. Our focus is on the role $C^1(M)$ eigenfunctions could play in the study of the Koopman operator as an operator on $L^2(\mu)$.

\begin{figure}[!t]
	\centering
	\begin{tikzpicture}[scale=0.55, transform shape]
\node [style={rect5}] (A1) at (\columnA, 1.5\rowA) {Assumption \ref{A:A1}};	
\node [style={rect5}] (Az) at (2\columnA, 1.5\rowA) {Assumption \ref{A:z}};	
\node [style={rect5}] (Agen) at (5\columnA, -4.5\rowA) {Assumption \ref{A:generating}};	
\node [style={rect5}] (Amltple) at (5\columnA, 1.5\rowA) {Assumption \ref{A:z_multple}};	
\node [style={rect2}] (1) at (0, 0) {Theorem \ref{thm:fol_1D_basic} -- Partitions from the leaves of a Koopman eigenfunction, and their conditional measures};
\node [style={rect2}] (2) at (2\columnA, 0) {Theorem \ref{thm:fol_1D_basemap} -- Return map and suspension flow over base leaf};
\node [style={rect2}] (3) at (4\columnA, 0) {Theorem \ref{thm:fol_1D_return} -- Ergodic properties of the return map.};
\node [style={rect2}] (4) at (\columnA, -1.5\rowA) {Theorem \ref{thm:fol_1D_flow} -- Restriction of the flow to the base leaf};
\node [style={rect2}] (5) at (3\columnA, -1.5\rowA) {Theorem \ref{thm:fol_1D_flow_eig} -- Restriction of the flow to the base leaf};
\node [style={rect6}] (6) at (4\columnA, -3\rowA) {Theorem \ref{thm:Main} -- Main result};
\node [style={rect7}] (7) at (2\columnA, -4.5\rowA) {Corollary \ref{corr:spectral_charac} -- Characterization of the quasiperiodic and chaotic components of the spectrum};
\node [style={rect7}] (8) at (4\columnA, -4.5\rowA) {Corollary \ref{corr:skew} -- Skew product structure of the dynamics, with toral rotations as drivers };
\node (9) at (5\columnA, -1.5\rowA) {};
\node (10) at (5\columnA, -3\rowA) {};
 \path [line2] (1) -- node [midway,below ] {} (2);
 \path [line2] (2) -- node [midway,below ] {} (3);
 \path [line2] (1) -- node [midway,below ] {} (4);	
 \path [line2] (2) -- node [midway,below ] {} (4);
 \path [line2] (3) -- node [midway,below ] {} (4);
 \path [line2] (4) -- node [midway,below ] {} (5);
 \path [line2] (4) -- node [midway,below ] {} (6);
 \path [line2] (5) -- node [midway,below ] {} (6);
 \path [line2] (6) -- node [midway,below ] {} (7);
 \path [line2] (6) -- node [midway,below ] {} (8);
 \path [line2, bend left] (3) -- node [midway,below ] {} (6);
 \path [line1] (A1) -- node [midway,below ] {} (1);
 \path [line1] (Az) -- node [midway,below ] {} (1);
 \path [line1] (A1) -- node [midway,below ] {} (2);
 \path [line1] (Az) -- node [midway,below ] {} (2);
 \path [line1] (A1) -- node [midway,below ] {} (3);
 \path [line1] (Az) -- node [midway,below ] {} (3);
 \path [line1] (Agen) -- node [midway,below ] {} (3);
 \path [line1] (A1) -- node [midway,below ] {} (4);
 \path [line1] (Az) -- node [midway,below ] {} (4);
 \path [line1] (Amltple) -- node [midway,below ] {} (Az);
 \path [line1] (Amltple) -- node [midway,below ] {} (9);
 \path [line1] (Amltple) -- node [midway,below ] {} (10);
 \path [line1] (9) -- node [midway,below ] {} (5);
 \path [line1] (10) -- node [midway,below ] {} (6);
 \path [line1] (Agen) -- node [midway,below ] {} (6);
\end{tikzpicture}
\vspace{1cm}
	\caption{Outline of the results. The main theoretical results of this article are presented here, along with the Assumptions they rely on. The dependence on the assumptions are indicated using thin grey arrows, while the logical dependencies of the main results are indicated y thick lack arrows. The Corollaries provide the most recognizable consequences of such eigenfunctions.}
	\label{fig:overview}
\end{figure}

\paragraph{Smooth eigenfunctions} So far we have seen the several utilities of studying the Koopman operator in the Hilbert space $L^2(\mu)$. The spectral theory for 1-parameter group of unitary operators is extremely powerful and complete in Hilbert spaces \cite{EisnerEtAl15, DGJ_compactV_2018, ValvaGiannakis2023cnstnt}. A drawback of this approach is that the resulting eigenfunctions are objects of the space $L^2(\mu)$. An $L^2(\mu)$ function is not a true function but an equivalence class of functions which are pairwise equal almost everywhere. Thus an $L^2(\mu)$-based approach would lack the usual interpretation of eigenfunctions as true functions. Each $L^2(\mu)$ equivalence class is represented by true functions. However, in most situations, these functions are known to be non-smooth. Topologically chaotic systems have a dense collection of periodic orbits. For such systems, $L^2(\mu)$-eigenfunctions with irrational eigenfrequencies are nowhere continuous. In spite of this obstacles, our focus is on $C^1(M)$-eigenfunctions. A dynamical system can have multiple or even uncountably many coexisting invariant regions, each supporting its own ergodic measure and displaying a different type of dynamics. A classical example is the standard Chirikov map \cite{Chirikov1971std, Chirikov1979std, DSSY_Mes_QuasiP_2016}. Note that $C^1(M)$ eigenfunctions are globally defined, and independent of a choice of ergodic measure.  Consider the following system
\begin{equation} \label{eqn:def:qpd}
	\begin{split}
		\frac{d \theta}{dt} &= \vec{\omega} \\
		\frac{d y}{dt} &= G(\theta, y)
	\end{split}
\end{equation}
called a \emph{quasiperiodically-driven} dynamical system (q.p.d.) \cite{DasMustAgar2023_qpd, DasEtAl2023traffic}. The $\theta$ above is a coordinate on a $d$-dimensional torus $\mathbb{T}^d$, and the $y$ on some manifold $\tilde{M}$. This results in a vector field on the product space $M = \mathbb{T}^d \times \tilde{M}$. The term $\vec{\omega}$ is a $d$-dimensional vector playing the role of \emph{angular speed}. For each $1 \leq j \leq d$ consider the function $z_j : M \to \cmplx$ given by $(\theta, y) \mapsto e^{\iota \theta_j}$. This makes $z_j$ a $C^1$ eigenfunction with eigenfrequency $\vec{\omega}_j$. The existence of these $d$ eigenfunctions do not depend on the choice of an invariant measure for the system above. In fact for a general flow $\Phi^t : M \to M$, any $C^1(M)$ eigenfunction is automatically an $L^2(\mu)$ eigenfunction for any choice of invariant measure $\mu$. As opposed to topologically chaotic systems, \eqref{eqn:def:qpd} is an explicit prototypical example of a dynamical system having $C^1$ eigenfunctions. One of our main results is that under very general assumptions, this is the only possibility : 

\begin{corollary} [Torus rotation driven dynamics] \label{corr:skew}
	Let Assumptions \ref{A:A1} and \ref{A:z_multple} hold. Then under a change of variables, the flow takes the form of a quasiperiodically driven dynamics \eqref{eqn:def:qpd}. 
\end{corollary}

Assumption \ref{A:A1} has already been stated. Assumption \ref{A:z_multple} appears later and is an assumption on the existence of $d$ $C^1(M)$ eigenfunctions along with a property called \emph{independence}. The eigenfrequencies of these eigenfunctions together form the vector $\omega$. Corollary \ref{corr:skew} is a direct consequence of \eqref{eqn:dof09} \blue{along with a special flow that we construct}. It is proved later in Section \ref{sec:conclus}. 

Note that the q.p.d. in \eqref{eqn:def:qpd} provides a direct factorization of the phase space $M$, and one of the components can be directly interpreted as the source of Koopman eigenfunctions. The statement of Corollary \ref{corr:skew} takes the opposite direction. It deduces the converse, from just the existence of $C^1(M)$ Koopman eigenfunctions. Under Assumptions \ref{A:A1} and \ref{A:z_multple}, any flow may be restructured as a skew-product dynamical system. The special class of skew product dynamical systems retain all the complexities of a general dynamical system but provide deeper insight into the influences of chaotic and quasiperiodic components. % \citep[e.g.][]{StepSkew1, StepSkew2, MelbourneStuart2011skew, DasJim17_chaos, DasMustAgar2023_qpd, DasEtAl2023traffic}.

If we consider  q.p.d. \eqref{eqn:def:qpd} once more, it is not apparent whether all the $C^1(M)$ or $L^2(\mu)$ eigenfunctions of this system arise from projections along the angular coordinates. There might be eigenfunctions of this system which depend on both the $\theta$ and $y$ coordinates. It turns out that under an additional general assumption, one can affirm that all the eigenfunctions are indeed only dependent on the angular coordinate $\theta$. In fact one obtains a complete characterization of the spectral components $\Disc$ and $\Disc^\bot$ from \eqref{eqn:L2_decomp} : 

\begin{corollary} [Discrete and continuous components] \label{corr:spectral_charac}
	Let Assumptions \ref{A:A1} and \ref{A:z_multple} hold, and consider the form \eqref{eqn:def:qpd} claimed by Corollary \ref{corr:skew}. If in addition Assumption \ref{A:generating} holds, then $\Disc$ and $\Disc^\bot$ respectively corresponds to the functions dependent only on the $\theta$ and $y$ variables respectively. 
\end{corollary}

Corollary \ref{corr:spectral_charac} is proved later in Section \ref{sec:conclus}. 
The Corollaries are a consequence of a series of results examining the deeper consequences of smooth Koopman eigenfunctions. See Figure \ref{fig:overview} for an overview of the main theoretical results and their logical dependencies. Theorems \ref{thm:fol_1D_basic}, \ref{thm:fol_1D_basemap}, \ref{thm:fol_1D_return} are presented in Section \ref{sec:partition}, where we explore the consequences of a single Koopman eigenfunction. In Section~\ref{sec:smooth_eig} we explore the additional consequences of the eigenfunctions being smooth. These consequences are stated precisely in Theorems \ref{thm:fol_1D_flow} and \ref{thm:fol_1D_flow_eig}. All these results are combined in Section \ref{sec:nested} to derive our main result Theorem \ref{thm:Main}. A collection of smooth, independent Koopman eigenfunctions are used to construct a series of vector fields, each of which is a successive spectral ``simplification'' of the original flow, and we end with a discrete map that isolates the continuous or chaotic component of the flow. The consequences of this deconstruction are explored in Section \ref{sec:conclus}. The proofs of all the results are presented in Section \ref{sec:proofs}.  

%-_-_-_-_-_-_-_-_-_-_-_-_-_-_-_-_-_-_-_-_-_-_-_-_-_-_-_-_-_-_-_-_-_-_-_-_-_-_-_-_-_-_-_-_-_-_-_-_-_-_-_-_-_-_-_-_-_-_-_-_-_-_-_-_-_-_-_-_-_-_-_-_-_-_-_-_-_-_-_-_-_-_-_-_-_-_
\section{Partition from Koopman eigenfunctions} \label{sec:partition}

Presently we shall restrict our attention to the usual $L^2(\mu)$ eigenfunctions, arising from the ergodic measure $\mu$ from Assumption \ref{A:A1}. We now describe a way to ``restrict'' the dynamics to the level sets of Koopman eigenfunction. Suppose we have a collection \{$z_{1}$, \ldots, $z_{d}$\} of $d$ Koopman eigenfunctions belonging to the class $C^r(M)$, $C(X)$ or $L^2(\mu)$. Let their eigenfrequencies  be $\omega_1,\ldots,\omega_d$ respectively. Set $\pi:X\to\TorusD{d}$ to be the map $x\mapsto \left(z_{1}(x), \ldots, z_{d}(x)\right)$. The fibres of $\pi$ create the following partition $\Fol$,
\begin{equation}\label{eqn:def:Fol}
\Fol:=\cup_{\theta\in\TorusD{d}}\Fol_\theta; \quad \Fol_\theta := \pi^{-1}(\theta),\ \forall\theta\in\TorusD{d}.
\end{equation}
If the $z_j$s are $L^2(\mu)$-eigenfunctions, then the partition $\Fol$ will only be a partition of the support $X$ of $\mu$. On the other hand, if the $z_j$ are $C(M)$ or $C(X)$-eigenfunctions, then \{$\Fol_\theta$ : $\theta\in\TorusD{d}$\} forms a partition of $M$ or $X$. In any of these cases, there is a family of conditional measures $\mu|\Fol_\theta$, each supported on the element $\Fol_\theta$ of this partition. Thus if $\Leb$ denotes the Lebesgue probability measure on $\TorusD{d}$, then
\[\theta\mapsto \int_{\Fol_\theta}\phi d(\mu|\Fol_\theta) \mbox{ is Borel-measurable } , \quad \int_M \phi d\mu = \int_{\TorusD{d}} \int_{\Fol_\theta}\phi d(\mu|\Fol_\theta) d \Leb(\theta) , \quad \forall \phi\in C(M)\]
When $d=1$, then $\pi$ coincides with a single Koopman eigenfunction $z$. We will state this as an assumption.

\begin{Assumption}\label{A:z}
$z:X\to S^1$ is a non-constant, $C(X)$ or $L^2(\mu)$ eigenfunction with eigenfrequency $\omega$. Let the partition created by its level sets be denoted as $\Fol$, similar to \eqref{eqn:def:Fol}. Set $\calL_1$ to be the leaf $\Fol_{0} \equiv z^{-1}(1)$.
\end{Assumption}

Our first result is about the conditional measures on the fibres of $z$. \blue{We shall make use of the measure theoretic notion of \emph{disintegration} of a measure with respect to (w.r.t.) a partition. It allows one to find the conditionals of a measure on measurable subsets of zero measure.}

\begin{theorem}\label{thm:fol_1D_basic}
Let Assumptions \ref{A:A1} and \ref{A:z} hold. Then for every $\theta\in [0,2\pi)$,
\begin{enumerate}[(i)]
	\item $\Phi^t$ maps the fibre $\Fol_\theta$ homeomorphically onto the fibre $\Fol_{\theta+\omega t \bmod 2\pi}$. 
	\item The images $\Phi^{2\pi t/\omega}(\Fol_\theta)$ are disjoint, for $0\leq t<1$, and their union is the whole set $X$.
	\item $\Phi^t_{*}\left(\mu|\Fol_{\theta}\right)$ = $\mu|\Fol_{\theta+\omega t \bmod 2\pi}$.
	\item The following map $\Iso$ is a continuous bijection.
	\begin{equation}\label{eqn:def:Iso1}	
	\Iso : \calL_1 \times [0,1)\to X; \quad \Iso(y,t) := \Phi^{t \cdot 2\pi/\omega} y
	\end{equation}
	Moreover $\Iso$ transforms the product measure $\left(\mu | \calL_1\right)\times \Leb$ into the invariant measure $\mu$. 
	\item If $z$ was an $L^2(\mu)$ eigenfunction instead, then these claims are true upto a set of $\mu$-measure $0$.
\end{enumerate}
\end{theorem} 

Theorem \ref{thm:fol_1D_basic} is proved in Section \ref{sec:proof:fol_1D_basic}. Note that Theorem \ref{thm:fol_1D_basic} is based on Assumption \ref{A:z}
which does not assume smooth eigenfunctions yet. \blue{Theorem \ref{thm:fol_1D_basic} demonstrates how an invariant measure disintegrates along the fibres of an eigenfunction, and the conditional measures along these fibres are preserved under the flow from fire to fire.} This is the first of a series of results that we derive on the consequences of Koopman eigenfunctions, to gain more insight on the phase spaces $X$ and $M$. 

\paragraph{Remark} $\Iso$ does not have a continuous inverse, other than the trivial case when $\Phi^{2\pi/\omega} \equiv Id$. To see this, let $a_0 \in \calL_1$ and $a_t$:= $\Phi^t(a_0)$ for every $t\in\real$. Then note that for every $t\in [0,T)$, $a_t$= $\Iso^{(1)}(a_0,t)$. Then, 
\[ \lim_{t\to 1^{-}} \Iso(a_0,t) = \lim_{t\to 1^{-}} \Phi^{t\cdot 2\pi/\omega}(a_0) = a_1 = \Phi^{2\pi/\omega}(a_0)\in \calL_1.\]
However note that $\left(\Iso\right)^{-1}(a_1)$ = $(a_1,0)$ is not the limit of $\left(\Iso\right)^{-1}(a_t)$ = $(a_0,t)$, as $t\to 1^{-}$.

\paragraph{Remark} The measure $\mu$ decomposes into conditional measures along the leaves of the foliation $\Fol_\theta$. Theorem \ref{thm:fol_1D_basic} thus says that these conditional measures can be also be obtained by pushing the conditional measure along any one leaf, under the flow. This observation will now be used to establish the isomorphism of the flow to a \emph{suspension flow}.

\begin{SCfigure}
\includegraphics[width=0.6\textwidth]{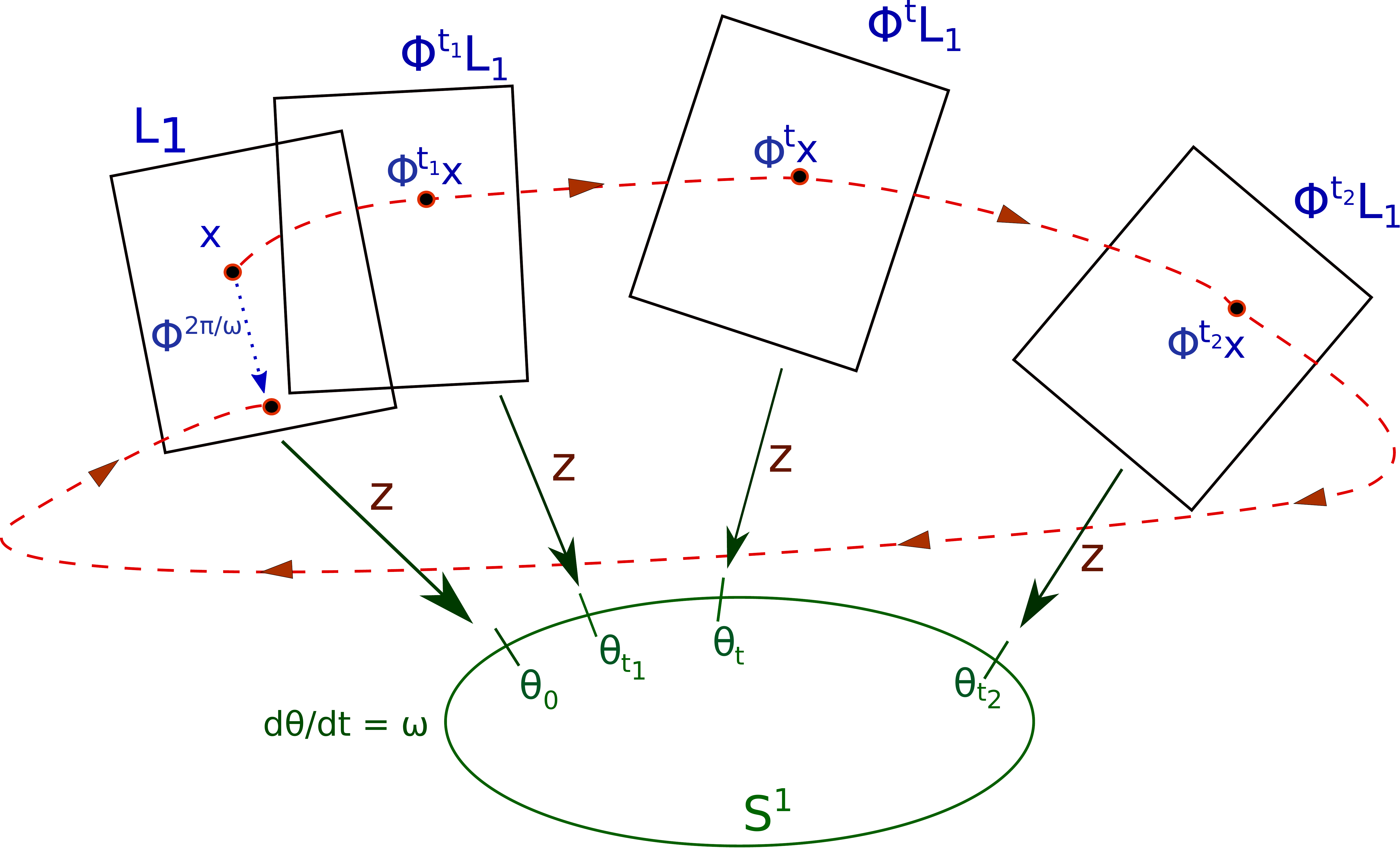}
\caption{Foliation $\Fol$ resulting from a Koopman eigenfunction $z:M\to S^1$ with eigenfrequency $\omega$. $\calL_1$ is a fixed leaf of the foliation as described in Assumption \ref{A:z} and \eqref{eqn:def:Fol}. It is assumed that $z|\calL_1 = 1$. \blue{ The eigenfunction $z$ takes values on the unit circle in the complex plane. This circle has been represented by its angular component $\theta$ on the circle. Each such $\theta$ is to be interpreted as the point $e^{\iota \theta}$ in $\cmplx$. } The flow $\Phi^t$ maps $\calL_1$ onto the entire leaf $\Phi^t(\calL_1)$, which corresponds to the inverse image $z^{-1} \paran{ \exp{\iota \theta_t} }$, where $\theta_t = t \omega \bmod 2\pi$. Thus $\Phi^{2\pi/\omega}$ maps $\calL_1$ into itself. One of the main results of Theorem \ref{thm:fol_1D_basic} is that $\Phi^t$ transforms the conditional measure $\mu|\calL_1$ into $\mu|\Phi^t(\calL_1)$.  }
\label{fig:fol_z}
\end{SCfigure}

Theorem \ref{thm:fol_1D_basic} and the remarks above suggest that the information of the dynamics can be expected to be contained in the base leaf $\calL_1$ of the eigenfunction. We next determine the exact nature of this information retention.

\paragraph{Suspension flows} Suspension flows are piecewise linear flows, and have been used extensively for constructing flows with various measure theoretic properties. We will begin by giving a restricted definition of a suspension flow which will be adequate for our purpose. Note that the map $\Return := \Phi^{2\pi /\omega}$ maps $\calL_1$ to $\calL_1$.  Let $\tilde{M}$ be the space $\calL_1\times \real / \sim$, where $\sim$ is an equivalence relation on $\calL_1 \times \real$ generated by the relations $(x,s+1)\sim (R_1(x),s)$. Then the suspension flow \emph{with base $\calL_1$, return map $\Return$ and height $1$} is a flow $\Susp^t$ on $\calL_1 \times \real / \sim$, defined as 
\begin{equation}\label{eqn:def:susp_flow}
\Susp^t(x,s) := \left( \Return^N(x), s+t-N \right), \quad \mbox{ where } N\in\integer \mbox{ is such that } 0\leq s+t-N<1.
\end{equation}
One can define, more generally, a suspension flow over some base space $\calL_1$, return map $\Return$, and height $h$ which is some function $h : \calL_1 \to \real^+$. Theorem \ref{thm:fol_1D_basemap} below shows that the original flow $\Phi^t$ is isomorphic to $\Susp^t$, and utilizes the special structure of $\Susp^t$ to make conclusions about the spectral properties of $\Phi^t$ in terms of the the return map $\Return$. 

\begin{theorem}\label{thm:fol_1D_basemap}
Let Assumptions \ref{A:A1} and \ref{A:z} hold, $\calL_1$ and $\Iso$ be as defined in \eqref{eqn:def:Iso1}, and $\Susp$ in \eqref{eqn:def:susp_flow}. Let $\nu$ the conditional measure of $\mu$ on $\calL_1$, obtained by its disintegration w.r.t. the foliation $\Fol$ from \ref{eqn:def:Fol}. Then 
\begin{enumerate}[(i)]
	\item $\Iso$ is a measure preserving, continuous conjugacy between the suspension flow $\Susp^t$ and the flow $\Phi^{2\pi t/\omega}$ :
     \[ \Iso\circ \Susp^t = \Phi^{2\pi t/\omega} \circ \Iso , \quad \forall t\in \real . \]
	\item Let $\zeta$ be an $L^2(\mu)$ (or $C(X)$ ) eigenfunction of $V$ with frequency $\omega'$, i.e., $V\zeta = \iota \omega' \zeta$. Then the map %$f:L \to L$ defined as
        \begin{equation} \label{eqn:def:base_map}
	   \tilde{\zeta} : \calL_1 \to \cmplx, \quad x \mapsto \zeta \paran{\Iso(x,0)}; \quad \forall x\in \calL_1,
	\end{equation}
	is an $L^2(\mu | \calL_1)$ (or $C(\calL_1)$ ) eigenfunction of the return map $\Return$ with eigenvalue  $e^{\iota 2\pi \omega'/\omega}$.
	\item Conversely, for any $L^2(\nu)$ (or $C(X)$ ) eigenfunction $\tilde{\zeta}$ of $\Return$ with eigenvalue $e^{\iota 2\pi \omega'/\omega}$, the map %$F:X\to\cmplx$ defined as
    \begin{equation}\label{eqn:def:tensor_with_base_map}
	\zeta :X\to\cmplx, \quad	\zeta\paran{ \Iso(y,s) } := e^{ \iota 2 \pi \omega' s/\omega} \zeta(y); \quad \forall y\in \calL_1, s\in [0,1)
    \end{equation}
    is an $L^2(\mu)$ (or $C(M)$ ) eigenfunction for $V$ with eigenfrequency $\omega'$.
\end{enumerate}
\end{theorem}

Theorem~\ref{thm:fol_1D_basemap} is proved in Section \ref{sec:proof:fol_1D_basemap}. It is a continuation of the analysis started in Theorem \ref{thm:fol_1D_basic}, by revealing the effect on the discrete spectrum. In broad terms, it states that the return map $\Return$ to the codimension-1 leaf $\calL_1$ has one less eigenfunction than the flow $\Phi^t$. We next find out that other spectral properties of $\Phi^t$ are retained in $\Return$. This requires a closer look at the components of the splitting in \eqref{eqn:L2_decomp}.

%-_-_-_-_-_-_-_-_-_-_-_-_-_-_-_-_-_-_-_-_-_-_-_-_-_-_-_-_-_-_-_-_-_-_-_-_-_-_-_-_-_-_-_-_-_-_-_-_-_-_-_-_-_-_-_-_-_-_-_-_-_-_-_-_-_-_-_-_-_-_-_-_-_-_-_-
\section{Multiple eigenfunctions} \label{sec:multiple}

We now discuss the algebraic structure intrinsic to the component $\Disc$ from \eqref{eqn:L2_decomp}

\paragraph{Generating frequencies} Given any two eigenfunctions $z_1, z_2$ with eigenfrequencies $\omega_1, \omega_2$, and any two integers $a,b$, the number $a\omega_1 + b\omega_2$ is also an eigenfrequency. This is because 
\[\begin{split}
    U^t \paran{ z_1^a z_2^b } &= \paran{ z_1^a z_2^b } \circ \Phi^t = \paran{ z_1^a \circ \Phi^t } \paran{ z_2^b \circ \Phi^t } = \paran{ z_1 \circ \Phi^t }^a \paran{ z_2 \circ \Phi^t }^b = \paran{ U^t z_1}^a \paran{ U^t z_2}^b = e^{\iota a \omega_1} z_1^a e^{\iota b \omega_2} z_2^b \\
    &= e^{\iota(a \omega_1 + b\omega_2)} \paran{ z_1^a z_2^b } .
\end{split}\]
	Thus integer linear combinations of frequencies are again frequencies, and products of eigenfunctions are again eigenfunctions. This makes the eigenfrequencies a module over the ring of integers. In particular, if the system has at least one nonzero frequency, then it has all harmonics of that frequency and thus infinitely many frequencies. A collection of eigenfrequencies are said to be independent if no integer linear combination of them is an integer.
	If the system has two independent frequencies, then all its frequencies are together dense on the real line. 
 
A collection of frequencies will be called a \emph{basis} or \emph{generating} set of eigenfrequencies if they are independent and all frequencies of the system can be generated by taking integer linear combinations of frequencies from this set. There is no unique choice of a basis, but all bases will have the same dimension $d$, called the quasiperiodicity dimension $d$. In finite-dimensional manifolds, the number $d$ is usually observed to be finite \citep[e.g.][Sec 3]{DasGiannakis_delay_2019}. \blue{A set of eigenfunctions corresponding to a generating set of frequencies form a \emph{generating set of eigenfunctions}. Note that such a collection automatically excludes the constant functions, since they correspond to the eigenfrequency $0$. The zero eigenfrequency cannot be included in a minimal spanning set of eigenfrequencies.} Our next focus will be on the component $\Disc^\bot$ from \eqref{eqn:L2_decomp}.

\paragraph{Mixing} A function $f\in L^2(\mu)$ is said to be (strongly) \emph{mixing} w.r.t. the ergodic dynamics $\paran{\mu, \Phi^t}$ if 
\[ \lim_{t\to\infty} \bracketBig{U^t f,g}_{\mu} = 0 , \quad \forall g\in L^2(\mu). \]
The function $f$ is said to be \emph{weakly mixing} w.r.t. the ergodic dynamics $\paran{\mu, \Phi^t}$ if 
\[ \lim_{T\to \infty} \frac{1}{T} \int_{t=0}^{T} \bracketBig{U^t f,g}_{\mu} dt = 0 , \quad \forall g\in L^2(\mu) . \]
The ergodic dynamics $\paran{\mu, \Phi^t}$ will be called weakly / strongly mixing if $f$ is weakly / strongly mixing for every $f\in L^2(\mu)$. Similar definitions can be made for discrete time dynamics, which are described by a map instead of a flow. \blue{Usually if the context is clear, we drop either $\mu$ or $\Phi^t$ when making a claim about mixing.}

Let $\alpha:\real^+ \to \real^+$ be a positive function such that $\lim_{t\to\infty} \alpha(t) = 0^+$. Then $\alpha$ is said to be a rate of (strong) mixing for the flow if
\[ \abs{ \bracketBig{U^tf, g}_{\mu} } = \bigO{ \alpha(t) }, \quad \mbox{ as } t\to\infty, \quad \forall f,g\in L^2(\mu). \] 
One can similarly derive notions of rates of strong and weak mixing for discrete and continuous time maps. Note that rate of mixing so defined, is not unique, but merely an upper bound for the asymptotic rate of decay. Our next result assumes :

\begin{Assumption}\label{A:generating}
The quasiperiodicity dimension $\dim_{\Disc}$ is finite. In other words the set $\Disc$ from \eqref{eqn:L2_decomp} is generated by $\dim_{\Disc}$ independent eigenfunctions $z_1,\ldots,z_{\dim_{\Disc}}$, for some $\dim_{\Disc}\geq 1$.
\end{Assumption}

\begin{theorem} \label{thm:fol_1D_return}
Let Assumptions \ref{A:A1}, \ref{A:z} and \ref{A:generating} hold, and  $\Return : \calL_1 \to \calL_1$ be the return map from Theorem \ref{thm:fol_1D_basemap}. Then
\begin{enumerate}[(i)]
	\item The discrete subspace ($\Disc$ from  \ref{eqn:L2_decomp}) of $\Return : \calL_1 \to \calL_1$ is generated by $\dim_{\Disc}-1$ eigenfunctions.
	
	Henceforth assume that $\dim_{\Disc}=1$. In that case,
	\item $\Return :  \calL_1 \to \calL_1$ has no eigenfunctions and hence is weakly mixing.
	\item $U^t | \Disc^\bot$ is strongly / weakly mixing iff $\Return$ is strongly / weakly mixing, in which case, these two systems have the same rates of mixing.
\end{enumerate}
\end{theorem}

While Theorems \ref{thm:fol_1D_basic} and \ref{thm:fol_1D_basemap} present the topological and measure theoretic consequences of a Koopman eigenfunction, Theorem \ref{thm:fol_1D_return} presents the consequence in terms of the discrete spectrum. \blue{Overall, we see that if we confine the dynamics to any fibre of a Koopman eigenfunction, then we get a simplification of the topology, the invariant measure, as well as the spectrum. The simplification of the spectrum is achieved by a reduction by one of the quasiperiodicity dimension.} This our first step towards understanding the spectral simplification in Corollary \ref{corr:spectral_charac}. Theorem \ref{thm:fol_1D_return} is proved in Section \ref{sec:proof:BI}. 

\paragraph{Remark} In the suspension flow stated in Theorem \ref{thm:fol_1D_return}, the flow is linear in-between two successive returns to the base $\calL_1$, and intuitively, the mixing in the dynamics can be expected to be solely contributed by the base map $\Return$. Theorem \ref{thm:fol_1D_return} shows that this is indeed the case if all the eigenfunctions are generated by a single eigenfunction, for then the dynamics of the return map has no non-trivial eigenfunction. This property is known as \emph{weak mixing} and will be dealt with in more details later. %

Theorems \ref{thm:fol_1D_basic} -- \ref{thm:fol_1D_return} performs a ``spectral reduction'' of the dynamics by considering only one leaf $\calL_1$ and looking at the return map to it. In this way, $\Disc$ has one less generator. This return map is however, discrete, and to spectrally reduce it further, we cannot apply these theorems again. We will next discuss how the vector field $V$ itself can be changed so as to be tangential to the foliation $\Fol$, but in such a way that the resulting flow is a reduced version of $\Phi^t$. For this, we need $z$ to be at least $C^1$ smooth. Having explored the consequences of the leaves of $C(X)$ or $L^2(\mu)$ Koopman eigenfunctions, we shall explore the further consequences of the eigenfunctions being smooth. 
%This will be the focus of the next section.

%-_-_-_-_-_-_-_-_-_-_-_-_-_-_-_-_-_-_-_-_-_-_-_-_-_-_-_-_-_-_-_-_-_-_-_-_-_-_-_-_-_-_-_-_-_-_-_-_-_-_-_-_-_-_-_-_-_-_-_-_-_-_-_-_-_-_-_-_-_-_-_-_-_-_-_-
\section{Smooth eigenfunctions} \label{sec:smooth_eig}

\begin{SCfigure}
\includegraphics[width=0.4\textwidth]{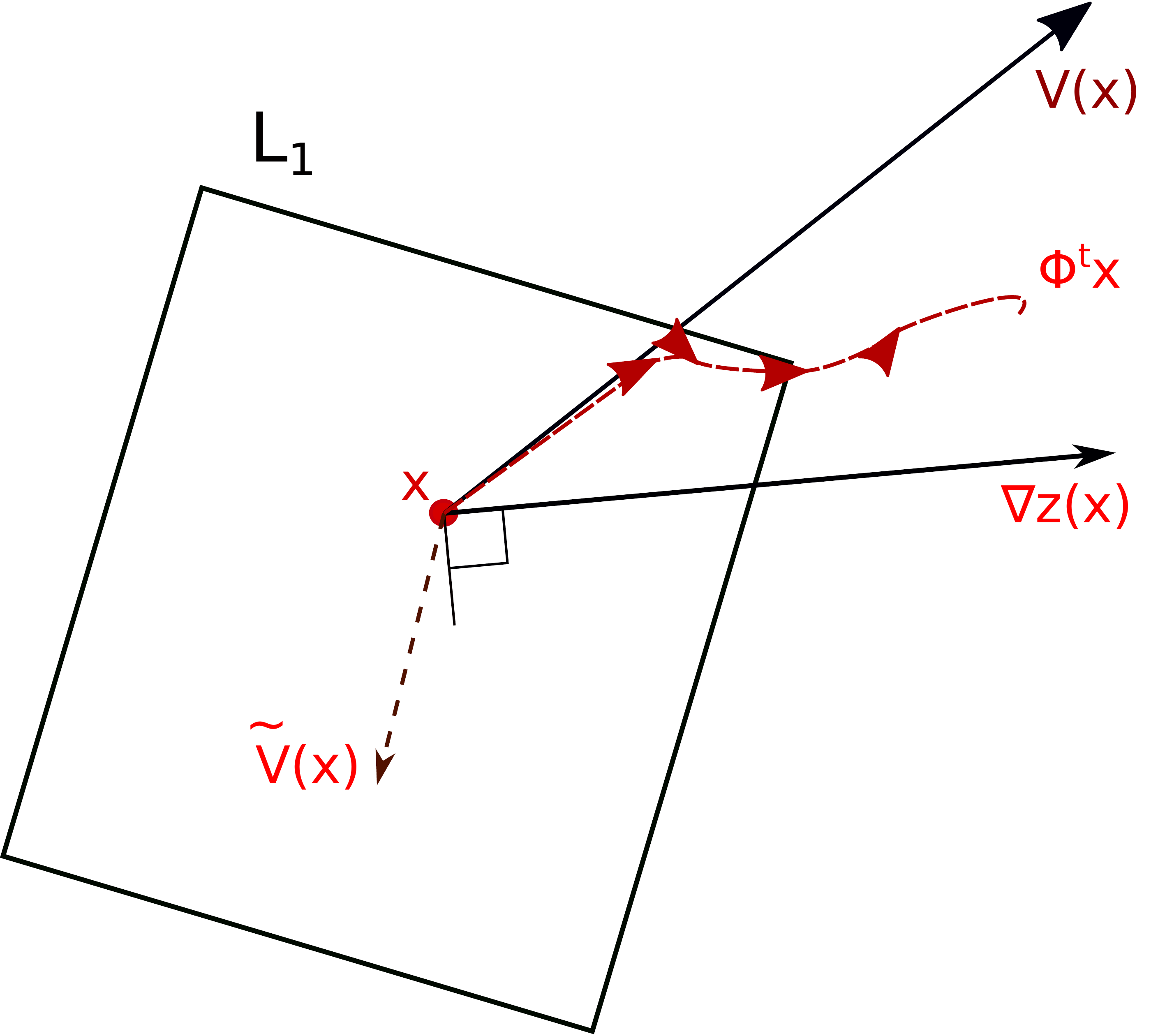}
\caption{Projection of the vector field $V$ along the foliation from $z$ results in a non-vanishing vector field $\tilde{V}$ \eqref{def:tildeV} tangential to the leaves of the foliation, see Theorem~\ref{thm:fol_1D_flow}. This projection is used in Theorem \ref{thm:Main} recursively to define the vector fields $V=V_1, V_2,\ldots, V_d$, and for $V_k$ has the property that it is tangential to the codimension $k-1$ foliation $\Fol^{(k-1)}$. See \eqref{eqn:def:Vk}  for precise definitions.}
\label{fig:def:Vk}
\end{SCfigure}

If we take a single non-constant $C^1(M)$ eigenfunction $z_1$, it gives a submersion of $M$ into $\TorusD{1}$ and this induces a codimension $1$ foliation. More generally, the following result shows that for $C^1(M)$ eigenfunctions,dynamic independence implies the differential property of submersion : 

\begin{lemma} \label{lem:submersion}
\citep[][Thm 6]{Das2023Lie} 
Let Assumption \ref{A:A1} hold \footnotemark[2], and $z_1,\ldots,z_m$ be $m$, independent, $C^1(M)$ eigenfunctions. Then the map $\pi:= (z_1, \ldots, z_m) :M\to\TorusD{m}$ is a submersion. 
\end{lemma}

\footnotetext[2]{The lemma's original statement in \cite{Das2023Lie} includes an assumption of an ergodic flow, which is not necessary.}

A $C^1$ Koopman eigenfunction $z:M\to S^1$ is also a submersion. Thus the leaves of the partition $\Fol$ are codimension-$1$ submanifolds. In Theorem~\ref{thm:fol_1D_flow} later, we will find a $C^1$ flow on each leaf of this foliation, and this flow retains all but one of the $d$ eigenfunctions. This process can be repeated $d$ times, and we will end up with a return map as in Theorem \ref{thm:fol_1D_return}, which retains the ergodic properties of the continuous spectrum of the flow. We will make this precise by describing a construction, that will decompose the flow in a manner similar to a suspension flow, as in Theorem \ref{thm:fol_1D_return}.  

To discuss the differential properties of $z$, we will need the concept of differential and gradient of an $S^1$ valued function $f:M\to S^1$. There is a way to define the differential $df$ as a $\real$-valued 1-form very similar to the case of a real valued function  \citep[see][]{Das2023Lie}. Moreover, for any choice of a Riemannian metric $g$ on $M$, the gradient to $f$, denoted as $\nabla_{\metric} f$, is the dual to $df$. Moreover, if one identifies $S^1$ with the unit circle on the complex place, then $f$ can be viewed as a map $f:M\to \cmplx$ and in that case,
\begin{equation}\label{def:grad_diff_Lie}
V(f) = i \bracket{V, \nabla_{\metric} f}_{\metric}  f = i df(V) f, \quad \forall  f\in C^1(M,S^1), \quad S^1 \subset \cmplx.
\end{equation}
We will next construct a vector field which is tangent to the foliation $\Fol$, and has one less independent eigenfunction than the original flow. 
\begin{theorem}\label{thm:fol_1D_flow}
Let Assumption \ref{A:A1} hold, and let $z:M\to S^1 \subset \cmplx$ be a $C^1$ eigenfunction with eigenfrequency $\omega$. Let $\Fol$ be the resulting foliation as in \eqref{eqn:def:Fol}, and $\calL_1$ be the leaf as in Assumption \ref{A:z}. Let $\nabla_\metric$ be the gradient of $z$ with respect to a metric $\metric$ on $M$, as described in \eqref{def:grad_diff_Lie}. Let $\tau$ be conformally transformed, so that $\nabla z$ is a unit vector field.  Define the vector field
\begin{equation}\label{def:tildeV}
\tilde{V} := V - \omega \nabla_{\metric} z
\end{equation}
Then the following hold. 
\begin{enumerate}[(i)]
    \item The vector field $\tilde{V}$ is tangent to the foliation $\Fol$. In particular, $\tilde{V}(z) = 0$ .
    \item $\tilde{V}$ is invariant under the flow $\Phi^t$, i.e., for every $x\in M$, every $t\in\real$, $\Phi^{t}_*|_x \tilde{V}(x)$ = $\tilde{V}\left( \Phi^{t} x \right)$. 
    \item The flow $\tilde{\Phi}^t$ := $\Phi^t_{\tilde{V}}$ generated by $\tilde{V}$ commutes with $\Phi^t$.
\end{enumerate}
\end{theorem}

\begin{SCfigure}
\includegraphics[width=0.6\textwidth]{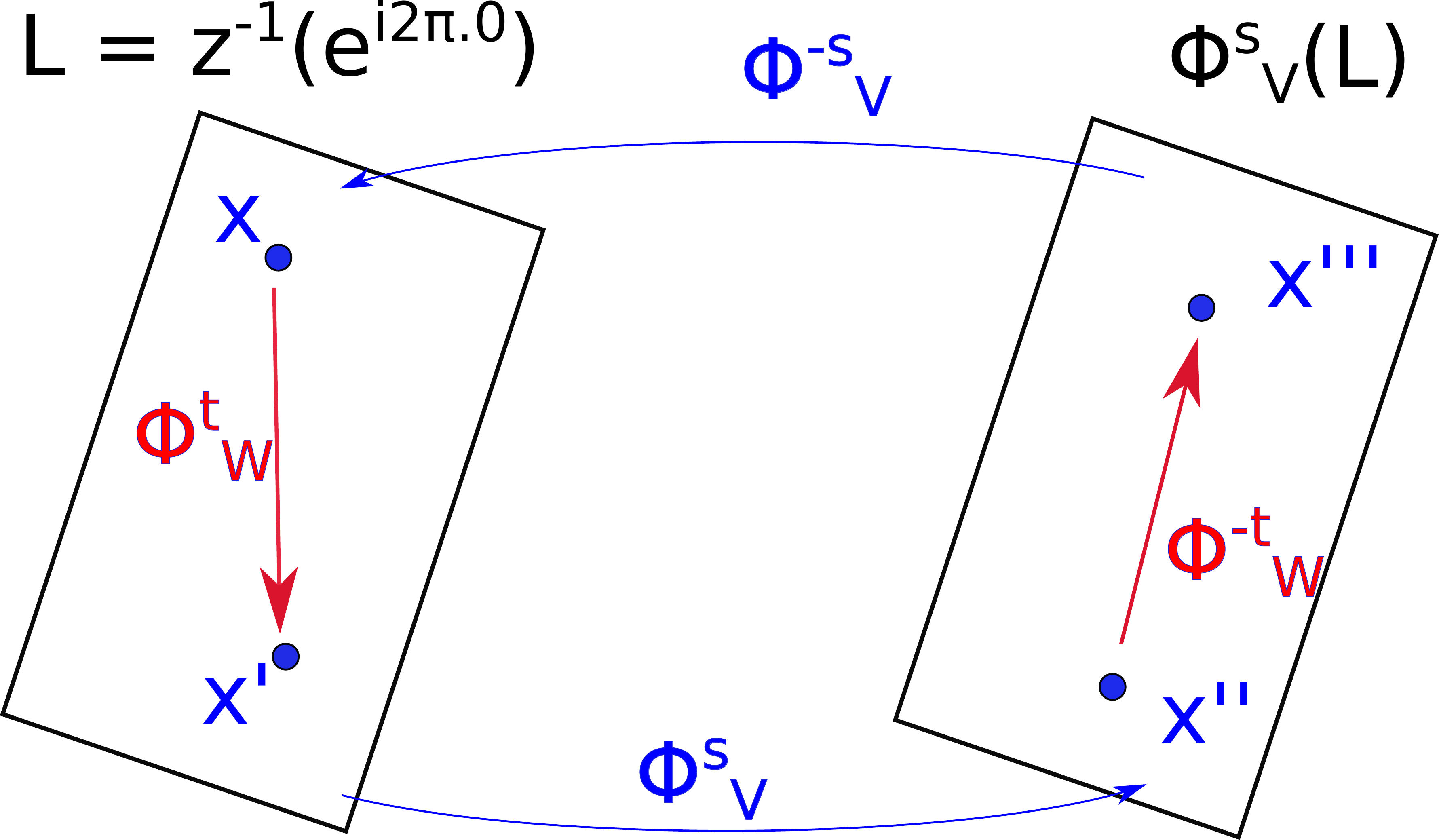}
%\;\begin{tikzcd}L\times [0,1) \arrow{r}{\Iso} &M \arrow{d}{W} \\T(L\times [0,1)) \arrow{r}{\Iso_*} &TM\end{tikzcd}
\caption{If $V,W$ are vector fields with corresponding flows $\Phi^t_V$ and $\Phi^t_W$, such that $\Phi^t_V$ preserves the vector field $W$, then the two flows commute. In Theorem~\ref{thm:fol_1D_flow}, $W$ corresponds to the vector field $\tilde{V}$, which is also tangent to the foliation resulting from the eigenfunction $z:M\to S^1$. %The left panel shows this commutation. The commutative diagram on the right exploits this commutation to build a skew-product representation of the flow $\Phi^t_V$, thus serving as a proof for Corollary~\ref{corr:skew}
}
\label{fig:comm_flow}
\end{SCfigure}

Theorem \ref{thm:fol_1D_flow} is proved in Section \ref{sec:proof:fol_1D_flow}. See Figure \ref{fig:comm_flow} for an illustration. zzz Recall the continuous bijection $\Iso: \calL \times  [0,1)\to M$ from \eqref{eqn:def:Iso1}. Our next Theorem \ref{thm:fol_1D_flow_eig} continues the analysis of the laminar flow $\tilde{\Phi}^t$. The following assumption will be used :

\begin{Assumption} \label{A:z_multple}
$z_1,\ldots,z_d$ are $d$ independent $C^1(M)$ eigenfunctions, with eigenfrequencies $\omega_1,\ldots,\omega_d$.
\end{Assumption}

\blue{ Note that the integer $d$ here is not necessarily the same as the quasiperiodicity dimension of the dynamics, which is the rang of the module of eigenfrequencies over the ring of rationals. This integer $d$ is definitely less than or equal to the quasiperiodicity dimension, and in our subsequent results we investigate the consequence of an equality. Recall that the eigenfrequencies inherit the module structure of the eigenfrequencies, but with pointwise multiplication playing the role of algebraic-addition. The quasiperiodicity dimension may also be interpreted in terms of eigenfunctions instead of eigenfrequencies. However, this number $d$ is different from the dimension of $\Disc$ as a vector space. This vector space is either one dimensional if there are no no-nontrivial eigenfunctions, or it is infinite. }

\begin{theorem} \label{thm:fol_1D_flow_eig}
    Under the same assumptions and notations of Theorem~\ref{thm:fol_1D_flow}, we further have
    \begin{enumerate} [(i)]
        \item Let $\tilde{\mu}$ be an invariant measure for  $\tilde{\Phi}^t$ restricted to $\calL$. Then for every $t\in [0,2\pi/\omega)$, $\Phi^t_* \tilde{\mu}$ will be an invariant measure for  $\tilde{\Phi}^t$ restricted to the leaf $ \Phi^t \paran{\calL}$.
        \item  Let $\tilde{\zeta}$ be an $L^2(\tilde{\mu})$ eigenfunction with eigenfrequency $\omega'$, for the flow $\tilde{\Phi}^t|_{\calL}$. Then its extension 
        \[ \zeta \paran{ \Iso(y,s) } := e^{\iota 2\pi s} \tilde{\zeta}(y), \quad \forall y\in \calL_1, \quad \forall s\in [0,1) . \]
        is an eigenfunction of $V$ with eigenfrequency $\omega'+\omega$.
        \item Suppose Assumption \ref{A:z_multple} holds. Then the metric $\metric$ may be chosen so that $\tilde{V}(z_j) = \iota \omega_j z_j$, for $j=1,\ldots,d$.
    \end{enumerate}
\end{theorem}

Theorem \ref{thm:fol_1D_flow_eig} is proved in Section \ref{sec:proof:fol_1D_flow_eig}. So far Theorems \ref{thm:fol_1D_basic},  \ref{thm:fol_1D_basemap} and \ref{thm:fol_1D_return} explore the consequences of a Koopman eigenfunction, the ergodic properties of the return maps, and how the mapping torus created by the return map recasts the flow as a suspension flow. Next Theorems  \ref{thm:fol_1D_flow} and \ref{thm:fol_1D_flow_eig} are based on the stronger assumption that the Koopman eigenfunctions are $C(M)$. In that case, the level sets of the eigenfunctions are codimension 1 submanifolds. presents a way to restrict the flow to the leaves, yet still retain many of the characteristics of the dynamics. In the next section we combine these results to construct a nested sequence of foliations and vector fields tangent to the foliations, each of which has quasiperiodicity dimensions successively $1$ less than the other.

%-_-_-_-_-_-_-_-_-_-_-_-_-_-_-_-_-_-_-_-_-_-_-_-_-_-_-_-_-_-_-_-_-_-_-_-_-_-_-_-_-_-_-_-_-_-_-_-_-_-_-_-_-_-_-_-_-_-_-_-_-_-_-_-_-_-_-_-_-_-_-_-_-_-_-_-
\section{Nested deconstruction} \label{sec:nested}

We henceforth assume Assumptions \ref{A:A1} and \ref{A:z_multple}. For every $k\in \braces{1, \ldots, d }$, let $1^{\times k}$ denote the vector $ (1,\ldots,1)$, which is also a point on $\TorusD{k}$. Set $\pi^{(k)} := \left( z_1, \ldots, z_k \right)$, which makes it a map $\pi^{(k)} : M \to \TorusD{k}$. By Lemma \ref{lem:submersion} $\pi^{(k)}$ is a submersion into $\TorusD{k}$. This lead to the foliations :
\[ \Fol^{(0)} := \{ X\}; \quad \Fol^{(k)}  := \cup_{\theta\in\TorusD{k}} \Fol^{(k)}_{\theta}; \quad \mbox{ where } \Fol^{(k)}_{\theta} := \left( \pi^{(k)} \right)^{-1}(\theta); \quad \forall \theta\in\TorusD{k}. \]
In fact $\Fol^{(k)}$ a a codimension-1 sub-foliation of $\Fol^{(k-1)}$. For each foliation we can define a base leaf
\[\calL_{0} := M; \quad \calL_k := \left( \pi^{(k)} \right)^{-1} \paran{1^{\times k}} .\]
Each $\calL_k$ is a codimension $k$ manifold and a leaf of the foliation $\Fol^{(k)}$. Note that the foliation $\Fol^{(k)}$ restricted to $\calL_{k-1}$ is also a codimension-1 sub-foliation, obtained from the submersion $z_k : \calL_{k-1}\to S^1$, and $\calL_k$ is the leaf corresponding to $z_k^{-1}(1)$. Similar to Theorem \ref{thm:fol_1D_basic}, one can inductively define conditional measures on the leaves $\calL_k$.
\begin{equation}\label{eqn:mu_k}
\mu_0:=  \mu; \quad \mu_k := \mu_{k-1}|\calL_{k}; \quad \forall k\in \braces{ 1, \ldots, d} .
\end{equation}

The gradient fields $\nabla_\metric z_1, \ldots, \nabla_\metric z_d$ are independent by Lemma \ref{lem:submersion}, thus a Riemannian metric $\tau$ can be chosen so that they are orthonormal to each other.  Now define $d$ vector fields $V_1,\ldots, V_d$, inductively as shown below.
\begin{equation}\label{eqn:def:Vk}
V_0 := V , \quad V_{k} := V_{k-1} - \langle V_{k-1}, \nabla_\tau z_k \rangle_{\tau} \nabla_\tau z_k, \quad \forall k\in \braces{ 1, \ldots, d}.
\end{equation}
The main result of this paper proves that the $V_k$ are tangential to these foliations $\Fol^{(k)}$, but at the same time, have $z_1,\ldots,z_k$ as eigenfunctions. It says even more.  Let $\Phi_k^t$ be the flow associated to $V_k$. Obviously, it is tangential to $\Fol^{(k)}$. Thus one has a return map 
\begin{equation} \label{eqn:def:Return_k}
    \Return_k := \Phi_{k}^{2\pi/\omega_k} : \calL_{k} \to \calL_{k}, \quad k\in \braces{1, \ldots, d} .
\end{equation}
\blue{We thus have a sequence of successively finer partitions $\Fol^{(k)}$, vector fields $V_k$ obtained by removing components of $V$ from the gradients of the first $k$ eigenfunctions, and the flows $\Phi^t_k$ induced by these vector fields. For each of the $k$-th level foliation, we have set apart the "first" leaf $\calL_k$. Our main result demonstrated how this nested construction also reflects a nested structure of the dynamics.}

\begin{theorem} [Main result] \label{thm:Main}
Let Assumptions \ref{A:A1} and \ref{A:z_multple} hold, and for every $k\in \braces{ 1, \ldots, d}$, $\pi^{(k)}$, $\Fol^{(k)}$ and $\calL_k$ be as defined above, $\mu_k$ be as in \eqref{eqn:mu_k},  let $V_k$ be the vector fields in \eqref{eqn:def:Vk}, $\Phi_k^t$ be the associated flow, and $\Return_k$ be the return maps from \eqref{eqn:def:Return_k}. Then for each $k$ in $\braces{1, \ldots, d}$,
\begin{enumerate}[(i)]
    \item $V_k$ is tangent to the foliation $\Fol^{(k)}$.
    \item $V_k z_j = 0$ if $1\leq j<k$,  and $V_k z_j = i \omega_j z_j$ if $k\leq j\leq d$, .
	\item The leaf $\calL_{k-1}$ is invariant under $\Phi_k$, i.e.,  $\Phi_k^t:\calL_{k-1} \to \calL_{k-1}$ is a $C^1$ flow. Moreover, $\Phi_k^{2\pi/\omega_k}:\calL_{k} \to \calL_{k}$.
	\item The map $\Iso^{(k)}$ defined below is a $C^1$ bijection.
	\begin{equation} \label{eqn:Td_translate}
		\Iso^{(k)}:\calL_k\times [0,1)^k \to M; \quad (x,t_1,\ldots,t_k) \mapsto \Phi_{1}^{2\pi t_1/ \omega_k} \cdots \Phi_{k}^{2\pi t_k/ \omega_k} x,
	\end{equation}
	\item The push-forward measure $\Iso^{(k)}_*(\mu_k \times \Leb^k)$ coincides with $\mu$.
	%\item The induced measure $\Phi_{k*}^{2\pi/\omega_k t} \mu_k$ on the leaf $\calL_k^t$ := $\Phi_k^{2\pi/\omega_k t} \calL_k$ coincides with its conditional measure w.r.t. the invariant measure $\mu$, namely, $\left(\mu|\Fol^{(k)}\right)(\calL_k^t)$.
	\item Suppose in addition, Assumption \ref{A:generating} holds with $\dim_{\Disc} = d$. Then $\Return_d: \calL_d \to \calL_d$ is weakly mixing. Moreover, $\Return_d$ has the same spectral measure as $U^t | \Disc^\bot$. In particular, $\Return_d$ is strongly mixing iff $U^t | \Disc^\bot$ is, in which case they have the same rates of mixing.
\end{enumerate}
\end{theorem}

\paragraph{Smoothness of eigenfunctions} Theorem~\ref{thm:Main} relies on the assumption of Koopman eigenfunctions being smooth. This condition is violated in several situations such as dynamics with SRB measures which have a dense collection of unstable periodic orbits. However any continuous time ergodic dynamical system may be spectrally approximated by smooth dynamical systems \cite{DGJ_compactV_2018, ValvaGiannakis2023cnstnt}. The spectral characteristic of a dynamical system is borne by the spectral measure of its associated Koopman operator $U$. If the dynamics is discrete time, then one can isolate the action of $U$ on $\Disc^\bot$. By a theorem of O. Knill \citep[][Thm 3.1]{Knill1997}, the linear dynamics $U \rvert_{\Disc^\bot}$ is conjugate to a translation on a topological group. Both these approximation results indicates that dynamical systems in which $\Disc$ is generated by a smooth set of eigenfunctions, is dense.

%Theorem~\ref{thm:Main} follows directly from the previous theorems. Claims~(i)--(iii) are consequences of Theorem~\ref{thm:fol_1D_flow} applied to the flow $\Phi^t_k$ on the manifold (leaf) $\calL_{k-1}$ and eigenfunction $z_k$, for $k\in \braces{ 1, \ldots, d}$. Claim~(iv) follows from Theorem~\ref{thm:fol_1D_basic}. Claim~(v) follows from Theorem~\ref{thm:fol_1D_return} and Theorem~\ref{thm:fol_1D_flow}~(iv). 

\paragraph{Proof of Theorem \ref{thm:Main}} 
We first prove Claims (i)-- (v). The proof will be by induction. We first check the base case corresponding to $k=1$. In this case, claims (i) and (iii) follow from Theorem \ref{thm:fol_1D_flow}~(i); claim (ii) follows from Theorem \ref{thm:fol_1D_flow_eig}~(iii); claim (iv) follows from Theorem \ref{thm:fol_1D_flow}~(iv); claim (v) follows Theorem \ref{thm:fol_1D_flow_eig}~(i). Next we verify the inductive step. So lets assume that the statement is true for every $k$ in the range $1\leq k \leq \bar{k}$ for some $1\leq \bar{k}<d$. Note that now $V_{\bar{k}}$ restricts to a vector field on $\calL_k$, which we can interpret as its own dynamical system. We can now construct a vector field $\tilde{V}_{\bar{k}}$ as in \eqref{def:tildeV}, which coincides with $V_{\bar{k}+1}$. The relation between $V_{\bar{k}+1}$ and $V_{\bar{k}}$ as listed in Theorems \ref{thm:fol_1D_flow} and \ref{thm:fol_1D_flow_eig} now complete the inductive step.

We next prove Claim (vi). Take any eigenvalue $\lambda$ of $\Return_d$. Again by unitarity, we have $\lambda = e^{\iota \omega'}$ for some eigenfrequency $\omega'$. Then by Theorem \ref{thm:fol_1D_basemap} and \eqref{eqn:def:tensor_with_base_map}, the flow $\Phi_d^t | \calL_{d-1}$ must have $\frac{\omega_{d} \omega'}{2\pi}$ as an eigenfrequency. By Theorem \ref{thm:fol_1D_flow_eig}~(ii), the flow $\Phi_{d-1}^t | \calL_{d-2}$ must have $\frac{\omega_{d} \omega'}{2\pi} + \omega_{d-1}$ as an eigenfrequency. Continuing this way, we get $\frac{\omega_{d} \omega'}{2\pi} + \omega_{d-1} + \cdots + \omega_1$ to be an eigenfrequency for the original flow. By Assumption \ref{A:generating}, there must be integers $a_1, \ldots, a_d$ such that 
\[\frac{\omega_{d} \omega'}{2\pi} + \omega_{d-1} + \cdots + \omega_1 = a_1 \omega_1 + \cdots + a_k \omega_d .\]
By the Independence of the generating frequencies, this is only possible if $\omega'/2\pi$ is an integer. But this would mean that  $\lambda=1$ a trivial eigenvalue for the return map. Thus the return map has no non-trivial eigenvalues and must be weakly mixing. The statements about weak or strong mixing, and their rates, now directly follow from Theorem \ref{thm:fol_1D_return}~(iii). This completes the proof of claim (v) and also of Theorem \ref{thm:Main}. \qed

\blue{This completes the proof and statement of our main results. Our assumptions have been very generic, except for the smoothness in Assumption \ref{A:z_multple}. We have seen how independent eigenfunctions lead to a succession of simplifications of the dynamics in more than one way.}
In the next section, we explore some consequences of Theorem~\ref{thm:Main}.

%-_-_-_-_-_-_-_-_-_-_-_-_-_-_-_-_-_-_-_-_-_-_-_-_-_-_-_-_-_-_-_-_-_-_-_-_-_-_-_-_-_-_-_-_-_-_-_-_-_-_-_-_-_-_-_-_-_-_-_-_-_-_-_-_-_-_-_-_-_-_-_-_-_-_-_-
\section{Conclusions} \label{sec:conclus}

One of the first important consequences of Theorem~\ref{thm:Main} is how the generating set of eigenfunctions interact with the change of variables $\Iso^{(d)}$ from \eqref{eqn:Td_translate}. For every $y\in \calL_d$ and $\paran{\theta_1, \ldots, \theta_d} \in [0,1)^d$, we have :
\begin{equation} \label{eqn:dof09}
    \begin{split}
        z_1 \circ \Iso^{(d)} \paran{ y, \theta_1, \ldots, \theta_d } &= e^{ \iota 2\pi \theta_1 } \\
        z_2 \circ \Iso^{(d)} \paran{ y, \theta_1, \ldots, \theta_d } &= e^{ \iota 2\pi \paran{\theta_1 + \theta_2} } \\
        \vdots &= \vdots \\
        z_d \circ \Iso^{(d)} \paran{ y, \theta_1, \ldots, \theta_d } &= e^{ \iota 2\pi \paran{\theta_1 + \cdots + \theta_d} } 
    \end{split} .
\end{equation}
These identities follow from Theorem \ref{thm:Main}~(ii). Equation \eqref{eqn:dof09} will be key to proving Corollary \ref{corr:skew}. % used to derive a simplified form for the dynamics, under a change of variables.

\paragraph{Proof of Corollary \ref{corr:skew}} Since $\Iso^{(d)}$ is a bijection, it has an inverse $\IsoInv$ whose components $\IsoInv_L : M \to \calL_d$ and  $\IsoInv_T : M\to [0,1)^d$ are defined uniquely by the relation below :
\begin{equation}\label{eqn:def:iso_inv}
	\IsoInv_L(z) = x, \; \IsoInv_T(x) = (t_1,\ldots,t_d) \quad \Leftrightarrow \quad \Iso^{(d)}(x,t_1,\ldots,t_d) = z,
\end{equation}
The following commutative diagram summarizes the relationship between $\Iso^{(d)}$, $\IsoInv$, $\IsoInv_L$ and $\IsoInv_T$.
\[\begin{tikzcd}[column sep = huge, row sep = large]
M \arrow{r}{\IsoInv_L} \arrow{d}[swap]{\IsoInv_T} \arrow[shift left=2]{dr}{\IsoInv} & \calL_d \\
\ [0,1) & \ [0,1)^d \times \calL_d \arrow{l}{\proj_T} \arrow{u}[swap]{\proj_L} \arrow[shift left=2]{ul}{\Iso^{(d)}}
\end{tikzcd}\]
Equipped with these notations, we can construct the flow $\Gamma^t$ indicated below.
\begin{equation} \label{eqn:def:Gamma_t}
    \begin{tikzcd} [row sep = large]
         M \arrow{rr}{ \Phi^t } & & M \arrow{d}[swap]{\Xi} \\
        \calL_d \times [0,1)^d \arrow{u}{ \Iso^{(d)} } \arrow[dashed]{rr}[swap]{\Gamma^t} & & \calL_d \times [0,1)^d
    \end{tikzcd}
\end{equation}
The dashed arrow indicates that the arrow is defined as the composition along the upper loop. This leads to an automatic satisfaction of the commutation. To verify that $\Gamma^t$ is indeed a flow, we have to verify that for any $t,s\in \real$, $\Gamma^{t+s} = \Gamma^t \circ \Gamma^s$. This follows from the commuting diagram below :
\[\begin{tikzcd} [row sep = large]
     M \arrow[dashed, bend left=30]{rrrr}{\Phi^{t+s}} \arrow{rr}{ \Phi^t } & & M \arrow[bend right=30]{d}[swap]{\Xi} \arrow{rr}{ \Phi^t } & & M \arrow{d}[swap]{\Xi} \\
    \calL_d \times [0,1)^d \arrow[dashed, bend right=30]{rrrr}[swap]{\Gamma^{t+s}} \arrow{u}{ \Iso^{(d)} } \arrow[dashed]{rr}[swap]{\Gamma^t} & & \calL_d \times [0,1)^d \arrow[bend right=30]{u}[swap]{ \Iso^{(d)} } \arrow[dashed]{rr}[swap]{\Gamma^s} & & \calL_d \times [0,1)^d
\end{tikzcd}\]
The two commuting squares on the left and right, as well as the square formed by the periphery, follow from \eqref{eqn:def:Gamma_t}. The middle loop created by the two opposite arrows hold since they represent inverses. The commutation in the upper loop holds because $\Phi^t$ is a flow. These five commutations come together to create the larger commutation diagram above. The consequence is the lower commutation loop, which is exactly the flow equation.

In spite of the vague nature of $\Phi^t$, the flow $\Gamma^t$ has a more explicit format : 
\begin{equation} \label{eqn:skew:1}
	\begin{split}
		\theta_1 &\mapsto \theta_1 + t\omega_1 \\
		\theta_2 &\mapsto \theta_1 + \theta_2 + t\omega_2 \\
		\vdots &\mapsto \vdots \\
		\theta_d &\mapsto \theta_1 + \cdots + \theta_d + t \omega_d \\
		y &\mapsto \paran{\Xi_L \circ \Phi^t \circ \Iso^{(d)}} \paran{ y, \theta_1, \ldots, \theta_d } 
	\end{split} ,
\end{equation}
for each $y\in \calL_d$ and $\paran{ \theta_1, \ldots, \theta_d } \in [0,1)^d $. Note that the change of variables 
\[ \paran{\theta_1, \ldots, \theta_d} \mapsto \paran{ \theta_1, \theta_1 + \theta_2, \ldots, \theta_1 + \cdots + \theta_d } , \]
converts \eqref{eqn:skew:1} into \eqref{eqn:def:qpd}, thus proving the claim of Corollary \ref{corr:skew}. \qed 

\paragraph{Proof of Corollary \ref{corr:spectral_charac}} Corollary \ref{corr:spectral_charac} relies on the fact that $\Iso^{(d)}$ is an invertible and measure preserving transformation between the two $\Phi^t$ and $\Gamma^t$. Thus there is a one-to-one correspondence between the discrete and continuous spectra of the two flows. According to the skew product nature of \eqref{eqn:skew:1}, the functions on $\calL_D \times [0,1)^d$ which only depend on the angular coordinates $[0,1)^k$ must lie in $\Disc$. Note that these functions have a set of $d$ generating functions. By Assumption \ref{A:generating} which we have included, this space is also maximal. Thus the space $\Disc$ is indeed the space $L^2 \paran{ [0,1)^d; \Leb }$. In that case its orthogonal component must be the space $L^2 \paran{ \calL_d; \mu_d }$. Equation \eqref{eqn:spectral_charac}, \eqref{eqn:Td_translate} describes how these functional spaces on $\calL_d\times [0,1)^d$ map back into $L^2(\mu)$. The decomposition $L^2(\mu) = \Disc \oplus \Disc^\bot$ from \eqref{eqn:L2_decomp} can now be explicitly written as 
\begin{equation} \label{eqn:spectral_charac}
	\Disc := \left\{ f\circ \IsoInv_T \;:\; f\in L^2 \paran{ [0,1)^d; \Leb } \right\} , \quad \Disc^\bot := \left\{ f\circ\IsoInv_L \;:\; f \in L^2 \paran{ \calL_d; \mu_d } \right\} .
\end{equation}
This completes the proof of the Corollary. \qed 

This completes the statement of all of our results. The rest of the paper will contain the proofs to the various theorems.

%-_-_-_-_-_-_-_-_-_-_-_-_-_-_-_-_-_-_-_-_-_-_-_-_-_-_-_-_-_-_-_-_-_-_-_-_-_-_-_-_-_-_-_-_-_-_-_-_-_-_-_-_-_-_-_-_-_-_-_-_-_-_-_-_-_-_-_-_-_-_-_-_-_-_-_-
\section{Proof of theorems and propositions} \label{sec:proofs}

%-_-_-_-_-_-_-_-_-_-_-_-_-_-_-_-_-_-_-_-_-_-_-_-_-_-_-_-_-_-_-_-_-_-_-_-_-_-_-_-_-_-_-_-_-_-_-_-_-_-_-_-_-_-_-_-_-_-_-_-_-_-_-_-_-_-_-_-_-_-_-_-_-_-_-_-
\subsection{Proof of Theorem \ref{thm:fol_1D_basic}} \label{sec:proof:fol_1D_basic}

Claim~(i) is a direct consequence of  \eqref{eqn:Def_koop_eigen}, and Claim~(ii) follows from Claim~(i). Let $\Borel(A)$ denote the Borel $\sigma$-algebras of of a topological space $A$.  Now define $T:= 2\pi/\omega$, so that
\[ \Phi^{sT} : \Fol_{\theta} \mapsto \Fol_{\theta+2s\pi \bmod 2\pi} , \quad \forall s\in\real, \quad \forall\theta\in S^1, \] 
and $X = \sqcup_{s\in [0,1)} \Fol_{2s\pi}$. For the proof of Claim~(iii), define a Borel measure $\MesL$ on $\calL_1$  as follows.
\[\MesL:= A\mapsto \int_0^1 \left(\mu|\Fol_{2s\pi}\right)\left(\Phi^{sT} A\right)ds, \quad \forall A\in \Borel \paran{ \calL_1 }.\]
In other words, given a measurable subset $A$ of the base leaf $\calL_1$, $\MesL(A)$ is the average of the conditional measures $\mu|\Fol_{2s\pi}$ on the images $\Phi^{sT} A \subseteq \Fol_{2s\pi}$. The first important observation is the following.
\begin{equation}\label{eqn:product_mes_tower}
\mu\left(\Iso\left(A\times [0,1)\right)\right) = \MesL(A), \quad \forall A\in \Borel(L). 
\end{equation}
To see why \eqref{eqn:product_mes_tower} holds, note that the set $S$ = $\Iso\left(A\times [0,1)\right)$ is the union $\cup_{s\in [0,1)} \Phi^{sT}(A)$. Thus by the definition of the conditional measures $\mu|\Fol_{2s\pi}$,  
\[\mu\left(S\right) = \int_X 1_Sd\mu = \int_0^1 ds\int_{\Fol_{2s\pi}} 1_S d(\mu|\Fol_{2s\pi}) = \int_0^1 \left(\mu|\Fol_{2s\pi}\right)\left(\Phi^{sT} A\right)ds = \MesL(A).\]
%This completes the proof of \eqref{eqn:product_mes_tower}. 
Now define $\tilde{\mu}$:=$\MesL\times\Leb$, the product measure on $\calL \times [0,1)$. We next prove that :
\begin{equation}\label{eqn:product_mes_push}
\Iso_{*}\tilde{\mu} = \mu. 
\end{equation}
To see why \eqref{eqn:product_mes_push} holds, first note that the following collection of rectangles form a sub-base for $\Borel(L\times [0,1))$.
%\begin{equation}\label{eqn:def:R}
\[ \{R_{A,i,n}:=A\times [\frac{i-1}{n},\frac{i}{n}) : \ A\in \Borel(L), \; n\in\num,\; 1\leq i \leq n \} \]
%\end{equation}
Therefore, it is sufficient to prove \eqref{eqn:product_mes_push} for each such rectangle $R= R_{A,i,n}$. Each of the $n$ rectangles $R_{A,j,n}$ for $j=1,\ldots,n$, have the same $\tilde{\mu}$ measure $\MesL(A)/n$. Moreover, 
 \[\Phi^{T/n} \left(\Iso(R_{A,j,n})\right) = \Iso(R_{A,j+1,n}),\quad \forall j=1,\ldots,n-1.\]
Therefore, by the flow-invariance of $\mu$, 
\[\mu\left(\Iso(R_{A,j,n})\right) = \mu\left(\Phi^{T/n}\left(\Iso(R_{A,j,n})\right)\right) = \mu\left(\Iso(R_{A,j+1,n})\right).\]
Therefore, all the images \{ $\Iso\left(R_{A,j,n}\right)$ : $j\in \braces{1,\ldots,n}$ \}  have the same $\mu$ measure. Therefore, 
%by \eqref{eqn:product_mes_tower}, 
\[\begin{gathered}
\mu(\Iso\left(R_{A,j,n}\right)) = \frac{1}{n} \sum_{k=1}^{n} \mu \left( \Iso\left(R_{A,k,n}\right) \right) = \frac{1}{n} \mu \left( \Iso\left( \sqcup_{k=1}^{n} R_{A,k,n}\right) \right) = \frac{1}{n} \mu\left(\Iso\left(A\times [0,1) \right)\right)\\
 \stackrel{\text{by \eqref{eqn:product_mes_tower}}}{=} \frac{1}{n}\MesL(A) =  \MesL(A)\times\Leb([\frac{j-1}{n},\frac{j}{n})) \stackrel{\text{by definition}}{=} \tilde{\mu}\left(R_{A,j,n}\right).
\end{gathered}\]
This completes the proof of \eqref{eqn:product_mes_push}. The proof of Claim~(iii) can now be completed. Fix an $A\subset L$ and set :
\[ \phi:[0,1)\to\real, \quad \phi(s) := \paran{ \mu|\Fol_{2s\pi} } \paran{ \Phi^{sT}(A) }, \quad \forall s\in[0,1). \]
The map $\phi$ is bounded and measurable. Next note that for every $[a,b)\subset [0,1)$,
\[\begin{split}
\frac{1}{b-a}\int_{a}^{b} \phi(t)dt &= \frac{1}{b-a} \int_{a}^{b} (\mu|\Fol_{2s\pi})(\Phi^{sT}(A)) ds = \frac{1}{b-a} \mu\left( \Iso(A\times[a,b]) \right) \\
&= \frac{1}{b-a}\tilde{\mu}\left(A\times [a,b]\right) \quad \mbox{[by \eqref{eqn:product_mes_push}]} \\
&= \MesL(A) \equiv \mbox{constant}.
\end{split}\]
Thus the integral on the left hand side (LHS) is independent of the integration limits $a$ and $b$, which is possible if only if $\phi$ is a constant function. Therefore, we have shown that
\begin{equation} \label{eqn:const_phi} 
    \phi(s) = (\mu|\Fol_{2s\pi})(\Phi^{sT}(A)) = \MesL(A) = \int_0^1 (\mu|\Fol_{2s\pi})(\Phi^{sT}(A)) ds, \quad \forall s\in [0,1).
\end{equation}
Using this identity we can further say that
\[\left( \Phi^{sT}_*(\mu | \calL_1) \right) (\Phi^{sT} A) = (\mu | \calL_1)(A) = \nu(A) =  (\mu|\Fol_{2s\pi})(\Phi^{sT} A), \quad \forall A\in\Borel(L).\]
The first equality follows from the definition of the push-forward of a measure. The last two equalities follow from \eqref{eqn:const_phi}. Since all sets in $\Borel(\Fol_{2s\pi})$ are of the form \{ $\Phi^{sT} A$ : $A\in\Borel(L)$ \}, the measures $\mu|\Fol_{2s\pi}$ and $\Phi^{sT}_*(\mu | \calL_1)$ must be the same and Claim~(iii) is proved. 

To prove Claim~(iv), first note that by \eqref{eqn:const_phi}, $\mu|\calL_1 = \nu$. Thus the measure $\tilde{\mu}$ equals $\mu|\calL_1 \times \Leb$. So by \eqref{eqn:product_mes_push}, $\Iso_{*} (\mu|\calL_1 \times \Leb) = \mu$, as claimed. The fact that $\Iso$ is a continuous bijection follows from Claims (i) and (ii). The details of Claim (v) will be omitted as they are exactly analogous. \qed

%-_-_-_-_-_-_-_-_-_-_-_-_-_-_-_-_-_-_-_-_-_-_-_-_-_-_-_-_-_-_-_-_-_-_-_-_-_-_-_-_-_-_-_-_-_-_-_-_-_-_-_-_-_-_-_-_-_-_-_-_-_-_-_-_-_-_-_-_-_-_-_-_-_-_-_-
\subsection{Proof of Theorem \ref{thm:fol_1D_basemap}} \label{sec:proof:fol_1D_basemap}

To prove Claim~(i), fix an arbitary point $(x,s)\in \calL_1 \times \real/\sim$. Without loss of generality (WLOG) we may assume that $0\leq s<1$. Then there is a unique integer $N$ such that $s+t-N\in [0,1)$. Then,
\[\begin{split} 
\Iso\circ \Susp^t (x,s) &= \Iso(\Return^N(x), s+t-N) = \Phi^{2\pi(s+t-N)/\omega} \Return^N(x) = \Phi^{2\pi(s+t-N)/\omega} \Phi^{2\pi N/\omega}(x)  \\
& = \Phi^{2\pi(s+t)/\omega} (x)=  \Phi^{2\pi t/\omega} \Phi^{2\pi s/\omega} (x) = \Phi^{2\pi t/\omega} \circ \Iso (x,s).
\end{split}\]
This confirms that $\Iso$ is indeed a conjugacy. The fact that it is a measure preserving conjugacy follows from Theorem~\ref{thm:fol_1D_basic}~(iv). This completes the proof of Claim (i).

To prove Claim~(ii), fix $x\in \calL_1$, and use the fact that $\Iso$ is a conjugacy (by Claim (i) ) to observe that 
\[ \Iso(\Return(y),0) = \Phi^{2\pi/\omega}\Iso(y,0), \quad \forall y\in \calL_1. %, \quad  \Iso(\Return(y),s) = \Phi^{2\pi s/\omega}\Iso(y,0), \; s\in [0,1). 
\]
%By assumption, for every $t\in\real$, $F(\Phi^t x) = e^{\iota \omega' t} F(x)$. 
Then by the definition of $\tilde{\zeta}$ in \eqref{eqn:def:base_map},
\[ \tilde{\zeta} \left( \Return(x) \right) =  \zeta\left( \Iso(\Return(x),0) \right) =  \zeta \left( \Phi^{2\pi/\omega} \Iso(x,0) \right) = e^{\iota 2\pi \omega'/\omega} \zeta\left( \Iso(x,0)\right) = e^{\iota 2\pi \omega'/\omega} \tilde{\zeta} (x), \]
%So $f$ is an eigenfunction of $\Return$ with eigenvalue $e^{\iota 2\pi \omega'/\omega}$, 
as claimed. To prove Claim~(iii), it has to be shown that for every $t>0$, every $x\in X$, $F\left( \Phi^{t} x \right)$ = $e^{\iota \omega' t} F(x)$. By Theorem~\ref{thm:fol_1D_basic}~(ii), $\Iso$ is a bijection, so there is $y\in \calL_1$ and $s\in [0,1)$ such that $x= \Iso(y,s)$. Let $n\in\integer$ such that $r=\frac{\omega t}{2\pi} + s -n \in [0,1)$. Then claim (iii) follows from the following substitutions.
\[\begin{split}
\zeta \left( \Phi^t x \right) &= \zeta \left( \Phi^t \Iso(y,s) \right) = \zeta \left( \Phi^t \Phi^{2\pi s/\omega} y \right) = \zeta \left( \Phi^{2\pi r/\omega} \Phi^{2\pi n/\omega} y \right) = \zeta \left(  \Phi^{2\pi r/\omega} \Return^{n} y \right) = \zeta \left( \Iso(\Return^{n}(y),r) \right)\\
& = e^{i 2\pi \omega' r/\omega} \tilde{\zeta} \paran{\Return^{n}(y)} = e^{i 2\pi \omega' r/\omega} e^{in2\pi \omega'/\omega} \tilde{\zeta}(y) = e^{ \iota 2 \pi \omega'(n+r)/\omega} \tilde{\zeta}(y) = e^{ \iota \omega' \paran{ t + s 2\pi/\omega } } \tilde{\zeta}(y) \\
& =e^{i \omega't} e^{ \iota 2 \pi \omega' s/\omega} \tilde{\zeta}(y) = e^{i \omega't} \zeta(\Iso(y,s)) = e^{i \omega't} \zeta(x). \\
\end{split}\]
This completes the proof of Theorem \ref{thm:fol_1D_basemap}. \qed

%-_-_-_-_-_-_-_-_-_-_-_-_-_-_-_-_-_-_-_-_-_-_-_-_-_-_-_-_-_-_-_-_-_-_-_-_-_-_-_-_-_-_-_-_-_-_-_-_-_-_-_-_-_-_-_-_-_-_-_-_-_-_-_-_-_-_-_-_-_-_-_-_-_-_-_-
\subsection{Proof of Theorem \ref{thm:fol_1D_return}} \label{sec:proof:BI}

We begin with a lemma that provides sufficient condition for a map to be strongly mixing.

\begin{lemma}\label{lem:mix_basis}
Let \{$\phi_k$ : $k\in\num$\} be a basis (not necessarily orthonormal) for $L^2(\mu)$, consisting of elements with uniformly bounded norm. Then 
\begin{enumerate} [(i)]
    \item the system is strongly  mixing iff for every $i,j\in\num$, $\langle U^T \phi_i, \phi_j \rangle$ $\to 0$ as $T\to\infty$. 
    \item If $\alpha:\real \to \real$ is a positive function such that $\lim_{|T|\to\infty} \alpha(T)$ = $0$ and for every $i,j\in\num$, $\left| \langle \phi_i, U^T \phi_j \rangle_\mu \right| = O(\alpha(T))$. Then $\alpha$ is a rate of mixing for $F$.
\end{enumerate}
\end{lemma}

Lemma \ref{lem:mix_basis} is proved in Section \ref{sec:proof:mix_basis}. Theorem \ref{thm:fol_1D_return} is now ready to be proved. Let us set $m := \dim_{\Disc}$. Take $m$ independent frequencies $\omega_1,\ldots,\omega_m$  of the flow $\Phi^t$, with $\omega=\omega_1$. Then by  Theorem \ref{thm:fol_1D_basemap}, $e^{\iota 2\pi \omega_2/\omega}, \ldots, e^{\iota 2\pi \omega_m/\omega}$ are $m-1$ independent eigenvalues of the base map $\Return$. Thus the discrete subspace of $\Return$ has at least $m-1$ independent eigenvalues. We will show that this is a maximal set. Consider any eigenvalue $\lambda$ of $\Return$, it must be of the form $e^{\iota 2\pi \omega'/\omega}$ for some $\omega'\in\real$. Then by Theorem \ref{thm:fol_1D_basemap}~(iii), $\omega'$ is an eigenfrequency of $\Phi^t$. Thus there are integers $a_1, \ldots, a_m$ such that $\omega' = \sum_{i=1}^{m} a_i \omega_i$. Therefore
\[ \lambda = \exp\left[ \frac{\iota 2\pi}{\omega} \sum_{i=1}^{m} a_i \omega_i \right] = \prod_{i=1}^{m} \left( e^{i a_i 2\pi\omega_i/\omega} \right)^{a_i}  = \prod_{i=2}^{m} \left( e^{i a_i 2\pi\omega_i/\omega} \right)^{a_i}. \]
The last equality follows from the fact that $\omega=\omega_1$. Thus $\lambda$ is generated by the $m-1$ eigenvalues, as claimed, proving Claim~(i). Claim (ii) is an immediate consequence of Claim~(i). 

By Claim~(i), the ergodic system $(\Phi^t,X,\mu)$ is measure theoretically isomorphic to a suspension flow $\Susp^t$ from \eqref{eqn:def:susp_flow}. Therefore these two systems have the same ergodic properties and it is sufficient to prove Claim~(iii) for the suspension flow. To that end, we use a special class $\Rect$ of functions, called \emph{rectangular functions}, instead of general functions in $L^2(L,\MesL\times \Leb)$. This class consists of functions of the form 
\[f(x,t) = 1_A(x)1_J(t),\]
where $1_A$, $1_J$ are characteristic functions of Borel measurable sets $A\subset X$ and $J\subset [0,1]$. We can further restrict $J$ to be an interval. These functions have a dense span in $L^2(L,\MesL\times \Leb)$ and by Lemma \ref{lem:mix_basis}, it suffices to check the mixing condition for such functions only. 

Let $f\in\Rect$. To simplify the analysis, we can assume without loss of generality that $f$ is of the form $A\times J$, with $J=[0,a)$ for some $a\in (0,1)$, and  $A\in\Borel(L)$.  Then,
\[ \Gamma^{n+s}(A\times J) = \left\{ \begin{array}{lcl} \Return^n(A)\times [s,s+a) &\mbox{ for } &0\leq s<1-a \\ \Return^n(A)\times [s,1) \cup \Return^{n+1}(A)\times[0,a+s-1) &\mbox{ for } &1-a \leq s<1 \end{array} \right. \quad \forall n\in\integer, \forall s\in [0,1). \]
 Let $g=1_B 1_K$ be another rectangular function. Let $t=n+s$ for some $n\in\integer$ and $s\in[0,1)$. Then
\[ \langle U^t f, g \rangle_\mu = \left\{ \begin{array}{lcl}
\MesL(B\cap\Return^n A)\Leb(K\cap [s,s+a)) &\mbox{ for } &0\leq s<1-a \\ 
\MesL(\Return^n(A)\cap B) \Leb([s,1)\cap K) + \MesL(\Return^{n+1}(A)\cap B)\Leb([0,a+s-1)\cap K) &\mbox{ for } &1-a \leq s<1 
\end{array} \right. \]
 Thus, we can conclude that
\begin{equation}\label{eqn:frhyt}
\frac{1}{\Leb(J)\Leb(K)} \langle U^t f, g \rangle_\mu \quad \mbox{is a convex sum of} \quad \MesL\left( \Return^{n}(A)\cap B \right), \MesL\left( \Return^{n+1}(A)\cap B \right).
\end{equation}
Both Claims (iii) and (iv) can now be seen to be true for rectangular functions, by virtue of \eqref{eqn:frhyt}. By Lemma \ref{lem:mix_basis}, they extend to all functions in $L^2(\MesL\times\Leb)$. 

We next prove Claim (ii). Take any two functions $f = \sum_{i\in\num} a_i \phi_i$  and $g = \sum_{i\in\num} b_i \phi_i$. Then
\[\begin{split}
    \abs{ \left\langle f, U^T g \right\rangle_{\mu} } & = \abs{ \left\langle \sum_{i\in\num} a_i \phi_i, U^T g \right\rangle_{\mu} } = \abs{ \left\langle \sum_{i\leq N} a_i \phi_i + \sum_{i>N} a_i \phi_i, U^T g \right\rangle_{\mu} } \\
    &\leq \abs{ \left\langle \sum_{i\leq N} a_i \phi_i , U^T g \right\rangle_{\mu} } + \abs{ \left\langle \sum_{i\geq N} a_i \phi_i , U^T g \right\rangle_{\mu} } \leq \abs{ \left\langle \sum_{i\leq N} a_i \phi_i , U^T g \right\rangle_{\mu} } + \norm{ \sum_{i\geq N} a_i \phi_i } \norm{U^T g} ,\\
    &\leq \abs{  \left\langle U^{-T} \sum_{i\leq N} a_i \phi_i , g \right\rangle_{\mu} } + \norm{ \sum_{i\geq N} a_i \phi_i } \norm{g} ,
\end{split}\]
where the last two steps utilizes the unitarity of $U$. Continuing, we get
\[\begin{split}
    \abs{ \left\langle f, U^T g \right\rangle_{\mu} } &\leq \abs{  \left\langle U^{-T} \sum_{i\leq N} a_i \phi_i , g \right\rangle_{\mu} } + \norm{ \sum_{i\geq N} a_i \phi_i } \norm{g} = \abs{  \left\langle  \sum_{i\leq N} a_i U^{-T} \phi_i , \sum_{j\in\num} b_j \phi_j \right\rangle_{\mu} } + \norm{ \sum_{i\geq N} a_i \phi_i } \norm{g} \\
    &= \abs{  \left\langle  \sum_{i\leq N} a_i U^{-T} \phi_i , \sum_{j\leq N} b_j \phi_j + \sum_{j\geq N} b_j \phi_j \right\rangle_{\mu} } + \norm{ \sum_{i\geq N} a_i \phi_i } \norm{g} \\
    &\leq \sum_{i,j\leq N} \abs{ a_i^* a_j  \left\langle U^{-T} \phi_i, \phi_j \right\rangle_{\mu} } + \abs{  \left\langle  \sum_{i\leq N} a_i U^{-T} \phi_i ,  \sum_{j\geq N} b_j \phi_j \right\rangle_{\mu} } + \norm{g} \norm{ \sum_{i\geq N} a_i \phi_i } \\
    &\leq \sum_{i,j\leq N} \abs{ a_i^* a_j  \left\langle U^{-T} \phi_i, \phi_j \right\rangle_{\mu} } + \norm{f} \norm{    \sum_{j\geq N} b_j \phi_j } + \norm{g} \norm{ \sum_{i\geq N} a_i \phi_i }
\end{split}\]
Fix an arbitrary $\epsilon>0$ and choose $N$ large enough so that the last two terms on the RHS above are less than $\epsilon$. For such $N$, the first term has a finite number $N^2$ of terms in the addition. Each of them are $\bigO{ \alpha(T) }$. Thus the entire RHS is $\bigO{ \alpha(T) } + 2\epsilon$. Since $\epsilon$ was arbitrary, the entire RHS and hence the LHS is $\bigO{ \alpha(T) }$. This completes the proof of claim (ii) and also of Theorem \ref{thm:fol_1D_return}. \qed

%-_-_-_-_-_-_-_-_-_-_-_-_-_-_-_-_-_-_-_-_-_-_-_-_-_-_-_-_-_-_-_-_-_-_-_-_-_-_-_-_-_-_-_-_-_-_-_-_-_-_-_-_-_-_-_-_-_-_-_-_-_-_-_-_-_-_-_-_-_-_-_-_-_-_-_-_-_-_-_-_-_-_-_-_-_-_
\subsection{Proof of Lemma \ref{lem:mix_basis} }  \label{sec:proof:mix_basis}

We begin with the proof of Claim (i). Note that by assumption, $C$:= sup\{ $\|\phi_j\|$ : $j\in\num$\} is finite. The ``only if'' part is obvious. So conversely, let for every $i,j\in\num$, $\langle U^T \phi_i, \phi_j \rangle$ $\to 0$ as $T\to\infty$. Let $f=\sum_{j\in\num} a_j \phi_j$. We first show that
\begin{equation} \label{eqn:pof03}
    \lim_{T\to\infty} \langle U^T \phi_i, f \rangle_\mu = 0, \quad \forall i\in\num .
\end{equation}
Let $\epsilon>0$. %It is enough to show that for $T$ large enough, $\left| \langle U^T \phi_i, f \rangle_\mu \right|<2\epsilon$. 
Choose $N\in\num$ large enough such that $\left\|\sum_{j=N+1}^{\infty} a_j \phi_j \right\| \leq \epsilon/C$. Then
\[\left| \langle U^T \phi_i, f \rangle_\mu \right| = \left| \sum_{j\in\num} \bar{a_j} \langle U^T \phi_i, \phi_j \rangle_\mu \right| 
\leq \left| \sum_{j=1}^{N} \bar{a_j} \langle U^T \phi_i, \phi_j \rangle_\mu \right| + \left| \sum_{j=N+1}^{\infty} a_j\langle U^T \phi_i, \phi_j \rangle_\mu \right|.\]
For $T$ large enough, by assumption on the $\phi_i$s, the first sum above becomes smaller than $\epsilon$. By the Cauchy Schwarz inequality and the fact that $U^T$ is a unitary operator,
\[  \left| \sum_{j=N+1}^{\infty} a_j\langle U^T \phi_i, \phi_j \rangle_\mu \right| = \left| \sum_{j=N+1}^{\infty} a_j\langle \phi_i, U^{-T} \phi_j \rangle_\mu \right| = \left| \langle \phi_i,U^{-T} \sum_{j=N+1}^{\infty} a_j  \phi_j \rangle_\mu \right| \leq \|\phi_i\| \left\|\sum_{j=N+1}^{\infty} a_j  \phi_j \right\| < \epsilon .\]
Therefore, the second sum is bounded by $\epsilon$. Thus for an arbitrary $\epsilon>0$, we have shown that $\left| \langle U^T \phi_i, f \rangle_\mu \right|<2\epsilon$. Thus \eqref{eqn:pof03} must be true.

Now let $g=\sum_{i\in\num} b_i \phi_i$. In a manner exactly analogous to the steps above, it can be shown that $\lim_{T\to\infty} \langle U^T g, f \rangle$= 0. The second claim of the lemma also follows in a similar manner. This completes the proof of Lemma \ref{lem:mix_basis}. \qed

%-_-_-_-_-_-_-_-_-_-_-_-_-_-_-_-_-_-_-_-_-_-_-_-_-_-_-_-_-_-_-_-_-_-_-_-_-_-_-_-_-_-_-_-_-_-_-_-_-_-_-_-_-_-_-_-_-_-_-_-_-_-_-_-_-_-_-_-_-_-_-_-_-_-_-_-_-_-_-_-_-_-_-_-_-_-_
\subsection{Proof of Theorem \ref{thm:fol_1D_flow}} \label{sec:proof:fol_1D_flow}

As shown in \citep[][Thm 2]{Das2023Lie} the gradient $\nabla_\metric z$ for such a submersion is orthogonal to the leaves of $\Fol$, with respect to (wrt) the $\metric$-metric. Thus Claim~(i) is equivalent to the identity $\langle \tilde{V}, \nabla_\metric z \rangle_\metric = 0$. By \eqref{def:grad_diff_Lie}, we have,
\[ \iota \omega z = V z = \iota \langle V , \nabla_\metric z \rangle_\metric z \quad \Rightarrow \quad \langle V , \nabla_\metric z \rangle_\metric = \omega. \]
Thus the inner product of $V$ with $\nabla_\metric z$ is constant everywhere in the $\metric$-metric. Therefore,
\[ \langle \tilde{V}, \nabla_\metric z \rangle_\metric = \langle V, \nabla_\metric z \rangle_\metric - \omega \langle \nabla_\metric z, \nabla_\metric z \rangle_\metric = \omega-\omega = 0 . \]
Thus $\tilde{V}$ is tangential to the leaves of $\Fol$. Since $z$ is constant on these leaves, $\tilde{V}(z) = 0$ everywhere, as claimed. 

We next prove Claim~(ii). Since the flow $\Phi^t$ preserves the leaves of the foliation, $\Phi^t_* TM \to TM$ preserves the sub-bundles $T\Fol$ and $T\Fol^\bot$. Thus there is a continuous function $\alpha:M\to\real$ such that 
\[\Phi^t_*|_x \nabla_\metric z(x) = \alpha(x) \nabla_\metric z(\Phi^t x), \quad \forall x\in M. \]
It will be shown that $\alpha(x)\equiv 1$. One of the fundamental facts from the theory of ordinary differential equations is that the flow $\Phi^t$ preserves the generating vector field $V$ :
\begin{equation} \label{eqn:Flow_commut_Gen}
     V(\Phi^t(x)) = d/ds|_{s=t} \Phi^s(x) = d/ds|_{s=0} \Phi^{t+s}(x) = \Phi^t_* (x)V(x) .
\end{equation}
By the definition of $\tilde{V}(x)$, $V(x) = \tilde{V}(x) + \omega \nabla_\metric z(x)$. Therefore,
\[ V(\Phi^t x) = \Phi^t_* (x)V(x) = \Phi^t_* (x)\tilde{V}(x) + \omega \alpha(x) \nabla_\metric z(\Phi^t x).   \]
Taking inner products on both sides with $\nabla_\metric z(\Phi^t x)$ gives
\[\omega = \langle V(\Phi^t x) , \nabla_\metric z(\Phi^t x) \rangle_\metric = \langle \Phi^t_* (x)\tilde{V}(x), \nabla_\metric z(\Phi^t x) \rangle_\metric + \omega \alpha(x) \langle\nabla_\metric z(\Phi^t x) , \nabla_\metric z(\Phi^t x) \rangle_\metric = 0 + \alpha\omega = \alpha\omega. \]
Thus $\alpha=1$ as claimed. This means that the gradient vector field $\nabla_\metric z$ is invariant under the flow, namely,
\[\Phi^t_*|_x \nabla_\metric z(x) = \nabla_\metric z(\Phi^t x), \quad \forall x\in M. \]
Therefore, since $\tilde{V}$ is a linear combination of the two flow invariant vector fields $V$ and $\nabla_\metric z$, $\tilde{V}$ is also invariant under $\Phi^t$. We summarize all these invariances below.
\begin{equation}\label{eqn:inv_vector}
\Phi^t_*|_x \tilde{V}(x) = \tilde{V}(\Phi^t x), \quad \Phi^t_*|_x V(x) = V(\Phi^t x), \quad \Phi^t_*|_x \nabla_\metric z(x) = \nabla_\metric z(\Phi^t x), \quad \forall x\in M. 
\end{equation}

We will next show that the two flows commute, i.e., $\tilde{\Phi}^t \Phi^s $ = $\Phi^s \tilde{\Phi}^{t}$ for every $s,t\in\real$. Let $x\in M$, and define $x'(s) := \tilde{\Phi}^s(x)$, $x''(s) := \Phi^t(x''(s))$ and $x'''(s) := \tilde{\Phi}^{-s}(x'(s))$. We have to show that $x'''(s) = \Phi^t(x)$ for every $s\in\real$. Since this is true for $s=0$, it is equivalent to prove that $\frac{d}{ds} x'''(s) = 0$. This follows from
\[\begin{split}
\frac{d}{ds} x'''(s) &= \frac{d}{ds} \tilde{\Phi}^{-s} \Phi^t \tilde{\Phi}^s(x) = \frac{d}{du}|_{u=s} \tilde{\Phi}^{-u} \Phi^t \tilde{\Phi}^s(x) + \frac{d}{dv}|_{v=s} \tilde{\Phi}^{-s} \Phi^t \tilde{\Phi}^v(x) \\
&= -\tilde{V}(x'''(s)) + \tilde{\Phi}^{-s}_*|_{x''(s)} \Phi^t_*|_{x'(s)} \tilde{V}(x'(s)) \\
& = -\tilde{V}(x'''(s)) + \tilde{\Phi}^{-s}_*|_{x''(s)} \tilde{V}(x''(s)), \quad \mbox{by } \eqref{eqn:inv_vector}, \\
& = -\tilde{V}(x'''(s)) + \tilde{V}(x'''(s)) = 0.
\end{split}\]
The last equality above again follows from the invariance of the vector field $\tilde{V}$ under the flow $\tilde{\Phi}^t$ that it generates. This completes the proof of Claim~(ii). This completes the proof of Theorem~\ref{thm:fol_1D_flow}. See Figure~\ref{fig:comm_flow} for an illustration. \qed

%-_-_-_-_-_-_-_-_-_-_-_-_-_-_-_-_-_-_-_-_-_-_-_-_-_-_-_-_-_-_-_-_-_-_-_-_-_-_-_-_-_-_-_-_-_-_-_-_-_-_-_-_-_-_-_-_-_-_-_-_-_-_-_-_-_-_-_-_-_-_-_-_-_-_-_-_-_-_-_-_-_-_-_-_-_-_
\subsection{Proof of Theorem \ref{thm:fol_1D_flow_eig}} \label{sec:proof:fol_1D_flow_eig}

The proof of Claim~(i) is analogous to the proof of Theorem~\ref{thm:fol_1D_basic}~(iii), and will be omitted. We next prove Claim~(iv). By \eqref{def:tildeV} we have
\[ V \zeta = \tilde{V} \zeta + \omega \nabla_\tau z \zeta . \]
The first differential operator acts only along the tangent subspaces of the leaves of $z$. The second differential operator acts along the vector field normal to the flow. Due to the local product structure,
\[ V \zeta = \iota \omega' \zeta + \iota \omega z \zeta . \]
Thus $\zeta$ satisfies the eigenfunction equation with eigenfrequency $\omega + \omega'$. If $\tilde\zeta$ is a $C^r$ eigenfunction, then so is $\zeta$. Suppose $\tilde\zeta$ is an $L^2(\tilde{\mu})$ eigenfunction. By Theorem \ref{thm:fol_1D_basic}~(iii) , since the flow $\Phi^t$ preserves the conditional measures of $\mu$ along the fibres of $\Fol$, the set $M'$ on which the eigenvalue equation holds has full $\mu$-measure. Thus, $\zeta$ is an $L^2(\mu)$ eigenfunction.

By Lemma \ref{lem:submersion}, the $d$ vector fields $\nabla_\metric z_1, \ldots, \nabla_\metric z_d$ are linearly independent everywhere on $M$, for any choice of a Riemannian metric $\metric$. Thus $\metric$ can be chosen so that these vector fields are mutually orthogonal and each of constant unit norm. As a result, using the identity \eqref{def:grad_diff_Lie} for the gradient, we get
\[ \tilde{V}(z_j) = V(z_j) - \omega \nabla_{\metric} z_1(z_j) = \iota \omega_j z_j - i \omega \langle \nabla_{\metric} z_1 , \nabla_{\metric} z_j \rangle_{\metric} z_j = \iota \omega_j z_j, \quad j=2,\ldots,d . \]
This completes the proof of Claim~(iii) and of Theorem~\ref{thm:fol_1D_flow_eig}. \qed

%-_-_-_-_-_-_-_-_-_-_-_-_-_-_-_-_-_-_-_-_-_-_-_-_-_-_-_-_-_-_-_-_-_-_-_-_-_-_-_-_-_-_-_-_-_-_-_-_-_-_-_-_-_-_-_-_-_-_-_-_-_-_-_-_-_-_-_-_-_-_-_-_-_-_-_-_-_-_-_-_-_-_-_-_-_-_
\bibliographystyle{\Path unsrt_inline_url}
\bibliography{\Path References,Ref}
\end{document}